\documentclass[preprint]{elsarticle}

\usepackage{amsmath,amssymb,latexsym,mathrsfs}
\usepackage{amsthm}
\usepackage{enumitem}
\usepackage{mathtools}
\usepackage{verbatim}
\usepackage{stmaryrd}
\usepackage{bussproofs}
\usepackage{tikz-cd}
\usepackage{natbib}
\usepackage{tikz}
\usepackage{xcolor}
\usepackage{adjustbox}
\usepackage[english]{babel}

\newtheorem{thrm}{Theorem}
\newtheorem{prop}{Proposition}

\newtheorem{lem}{Lemma}

\theoremstyle{definition}
\newtheorem{defn}{Definition}[section]

\theoremstyle{remark}
\newtheorem{exmp}{Example}
\newtheorem{rem}{Remark}

\newcommand{\m}{\mathcal}
\newcommand{\s}{\mathscr}
\newcommand{\ca}{\mathcal}

\DeclareMathOperator{\Th}{Th}

\DeclareMathOperator{\Diag}{Diag}

\DeclareMathOperator{\VS}{VS}
\DeclareMathOperator{\fin}{fin}

\DeclareMathOperator{\Mod}{Mod}

\DeclareMathOperator{\Ob}{Ob}
\DeclareMathOperator{\Mor}{Mor}
\DeclareMathOperator{\Alg}{Alg}

\DeclareMathOperator{\Log}{Log}
\DeclareMathOperator{\A}{A}

\newcommand{\Set}{\mathbf{Set}}

\newcommand{\Pos}{\mathbf{Pos}}

\newcommand{\Cat}{\mathbf{Cat}}

\newcommand{\FP}{\mathbf{FP}}

\newcommand{\FA}{\mathbf{FA}}

\DeclareMathOperator{\Sg}{Sg}

\newcommand{\db}[1]{\llbracket #1 \rrbracket}
\newcommand{\blank}{\rule{.2 cm}{0.15mm}}

\newcommand{\con}{\otimes}
\newcommand{\bigcon}{\bigotimes}
\newcommand{\propo}{\mathrm{prop}}

\begin{document}
	
	\begin{frontmatter}
	\title{Fibered universal algebra for first-order logics \tnoteref{t1}}
	\author[a]{Colin Bloomfield\corref{cor1}}
	\ead{colinbloomfield1@gmail.com}
	\author[b]{Yoshihiro Maruyama}
	\ead{yoshihiro.maruyama@anu.edu.au}	
	\cortext[cor1]{Corresponding author}  
	\address[a]{Department of Mathematics, Vanderbilt University, 1326 Stevenson Center Ln, Nashville, TN 37212, USA}
	\address[b]{School of Computing, Australian National University, 108 North Rd, Acton ACT 2601, Australia}
	\begin{abstract}
		We extend Lawvere-Pitts prop-categories (aka.\ hyperdoctrines) to develop a general framework for providing ``algebraic'' semantics for nonclassical first-order logics. This framework includes a natural notion of substitution, which allows first-order logics to be considered as structural closure operators just as propositional logics are in abstract algebraic logic. We then establish an extension of the homomorphism theorem from universal algebra for generalized prop-categories and characterize two natural closure operators on the prop-categorical semantics. The first closes a class of structures (which are interpreted as morphisms of prop-categories) under the satisfaction of their common first-order theory and the second closes a class of prop-categories under their associated first-order consequence. It turns out, these closure operators have characterizations that closely mirror Birkhoff's characterization of the closure of a class of algebras under the satisfaction of their common equational theory and Blok and Jónsson's characterization of closure under equational consequence, respectively. These ``algebraic'' characterizations of the first-order closure operators are unique to the prop-categorical semantics. They do not have analogs, for example, in the Tarskian semantics for classical first-order logic. The prop-categories we consider are much more general than traditional intuitionistic prop-categories or triposes (i.e., topos representing indexed partially ordered sets). Nonetheless, to our knowledge, our results are still new, even when restricted to these special classes of prop-categories.
	\end{abstract}
	\end{frontmatter}
	
	\section{Introduction}
	
	There is a vast and ever-growing expanse of nonclassical logics and with it, a need to develop a general framework to classify, prove metalogical results, and develop semantics for \textit{classes} of logics. For propositional logics, the field of abstract algebraic logic (AAL) provides such a framework rooted in universal algebra. However, for first-order logics, there is not a clear consensus on what the unifying algebraic perspective should be \cite[p.\ 104]{Font2016}. Two main challenges mentioned in the literature are (1) a lack of \textit{general} algebraic approaches to providing semantics for first-order logics and (2) a lack of a notion of formula substitution which preserves logical consequence \cite[p.\ 106]{Font2016}.
	
	To address (1), we use Lawvere's theory of hyperdoctrines \cite{Lawvere1969} to provide ``algebraic'' semantics for nonclassical first-order logics. A hyperdoctrine consists of a ``base'' category $ \m{C} $ with finite products whose objects are called types, and a contravariant functor $ P \colon \m{C} \to \Cat $ which associates a category of ``attributes'' to each type. A key insight of the hyperdoctrine approach to modeling first-order logics is to consider terms and formulas \textit{in context}, that is, with a list of variable/sort pairs that at least include those occurring in the expression. Then in the hyperdoctrine representation of first-order syntax, contexts are types, (lists of) terms are morphisms between types, the attributes associated to a context are the formulas in that context and derivations of one formula from another in the same context are morphisms in their category of attributes. We note the use of contexts is ubiquitous in classical first-order model theory where it allows one to unambiguously interpret terms as term functions and formulas as definable sets. However, it is scarcely formally presented there, which is unfortunate because contexts allow one to give an elegant point-free/valuation-free definition of satisfaction.\footnote{Lawvere's point-free semantics also resolves a philosophical issue with assigning objects to variables discussed in \cite[p.\ 10]{Walsh2018}.}

	Lawvere is also known for the view that theories are categories, possibly with some additional properties \cite{Lawvere1963,Lawvere1969}; this perspective is
	taken, for example, in the influential textbook on categorical model
	theory by Makkai and Reyes \cite[p.\ 168]{Makkai1977} and in Johnstone’s tome on
	topos theory \cite[p.\ 841]{Johnstone2002}. The two approaches (i.e., hyperdoctrines and theories-as-categories) are essentially equivalent for
	intuitionistic logic and its relatives such as coherent logic \cite{Frey2015}. This paper nevertheless takes the hyperdoctrine approach, since the framework of the paper is expected to cover monoidal logics as well as those cartesian logics, such as first-order linear and other substructural logics, for which the hyperdoctrine approach works better than the theories-as-categories approach. At the same time it
	should be noted that the theories-as-categories approach is
	encompassed by the hyperdoctrine approach in the sense that type theories (e.g., higher-order intuitionistic logic) can be represented as base categories in hyperdoctrines (e.g., base toposes in subobject hyperdoctrines); in this particular sense the two approaches are not necessarily in conflict with each other.

	 Since we model truth and not proof, we use poset-valued instead of more general category-valued hyperdoctrines for our semantics called prop-categories by Pitts \cite{Pitts2000}. That is, a prop-category is a pair $ (\m{C},P) $, where $ \m{C} $ is a category with designated finite products and $ P \colon \m{C} \to \Pos $ is a contravariant functor to $ \Pos $, the sub-category of $ \Cat $, whose objects are partially ordered sets and morphisms are monotone maps. The prop-categories in \cite{Pitts2000} are assumed to satisfy additional adjointness conditions so that they interpret the quantifiers, equality predicate and connectives of intuitionistic predicate logic. To model a broad variety of logics, we allow arbitrarily many quantifier and operation symbols in a first-order language $ \s{L} $ and we replace the adjointness conditions with weaker algebraic ones. Despite this generalization, the semantics of the logics we consider can still be encoded in a $2$-category of prop-categories we denote $ \FA_{\s{L}} $ (or just $ \FA $), where roughly theories are prop-categories, structures are prop-category morphisms, and structure-preserving maps are $2$-cells.

	Moreover, the morphisms in $ \FA $ determine a natural action of formula substitution which also addresses issue (2). A standard reference on AAL \cite[p.\ 106]{Font2016} says the following about substitution in first-order logics in contrast with substitution in propositional logics:
	\begin{quote}
		\ldots the role of sentential variables is played by atomic formulas, and the concept of substitution itself becomes quite different. In first-order languages we replace an individual (free) variable by a term; we do not replace an atomic formula, inside a more complex one, by an arbitrary formula, which is what is done in a sentential language. As a consequence, no substitution-invariance (structurality) holds.
	\end{quote}
	However, if one unpacks the action of morphisms on theories it essentially\footnote{Technically, in the case of single-sorted classical predicate
	logic, for a fixed $n \in \mathbb{N}$, we replace each $m$-ary relation symbol with an arbitrary formula-in-context, whose context has $nm$ variables. This is a natural generalization of the propositional case, since propositional atomic formulas are also nullary relation symbols.} does that. A restricted form of this notion of formula substitution appears in \cite[p.\ 191]{Church1956}, though it seems to be missing from many contemporary textbooks on mathematical logic. This is not surprising, since the account in \cite{Church1956} for classical first-order logic is quite complicated without the help of contexts.

	Much like the study of algebraic semantics for first-order logic, currently, the study of prop-categorical semantics for logics (in contrast with type theories) appears to be focused on studying particular logics and not on establishing a general framework in the spirit of AAL. An early exception is the work of Hiroyuki Shirasu, who in \cite{Shirasu1995} provides complete prop-categorical semantics for substructural predicate logics, and uses this semantics to prove their disjunctive and existence properties. Also, in \cite{Shirasu1998} Shirasu develops both prop-categorical and a general metaframe semantics for first-order modal logics and proves they are dual to each other.  More recently, in \cite{Maruyama2021}, complete hyperdoctrine semantics for substructural predicate logics is also given and used to present a unified account of logical translations, including Kolmogorov’s double negation translation of classical logic into intuitionistic logic and Girard’s exponential translation of intuitionistic logic into linear logic. In contrast to these existing works on general prop-categorical semantics we model nonclassical quantifiers and focus on (1) establishing minimal conditions on prop-categories and first-order logics so that structures are morphisms in a $ 2 $-category of prop-categories and (2) developing general algebraic/categorical results on the resulting prop-categorical semantics.

	Towards (1),  we define a weak logic $ \m{L}^{m} $ such that whenever a logic $ \m{L} $ is stronger than $ \m{L}^{m} $, each $ \m{L} $-theory $ T $ defines a classifying prop-category $ (\m{C}_{T},P_{T}) $ in $ \FA$, which contains a generic $ T $-model $ G $ in $ (\m{C}_{T},P_{T}) $. Moreover, $ \m{L}^{m} $ has a complete semantics with respect to the class of all of its classifying prop-categories (Theorem \ref{thm: FAsL is a complete semnatics for Lm.}), as do many natural extensions of $ \m{L}^m $ by the addition of various sequent rules (Theorem \ref{thm: General Completeness Theorem}). We then develop an algebraic view of the prop-categorical semantics of a logic $ \m{L} $ in analogy with the algebraic approach to the semantics of propositional logics. In this view, a full sub-2-category $ \FA_{\m{L}} $ of $ \FA $ replaces the usual quasivariety of algebras forming the algebraic semantics of a propositional logic and we define the kernel of a morphism (Definition \ref{defn: Kernel of a morphism.}) to give an algebraic/morphic view of logical entailment. These results are directly applied to formulate and prove the later ``fibered'' universal algebraic results.

	After singling out the class of first-order logics which have nice categorical semantics, in the second part of the paper, we develop a ``fibered universal algebra'' for $\FA $. The internal logic of a prop-category $ (\m{C},P)$ in $ \FA $ is defined and shown to be an $ \m{L}^{m} $-theory $ T $, such that $ (\m{C}_{T},P_{T}) \equiv (\m{C},P) $ in  $ \FA $. Subprop-categories and morphic images of prop-categories are defined and $ \FA $ is shown to have arbitrary products. The logical properties of products and subprop-categories are considered which mirror the logical properties of submatrices and products of matrices respectively in the matrix semantics of AAL \cite[p. 183]{Font2016}.
	
	We then prove a ``fibered'' analogue of the homomorphism theorem from universal algebra in two parts (Theorem \ref{thm: Fibered homomorphism theorem}). The first part, says that each morphism $ F \colon (\m{C},P) \to (\m{D},Q) $ in $ \FA_{\m{L}} $ factors through a ``quotient'' prop-category which is also in $ \FA_{\m{L}} $. The second part characterizes when, given morphisms $ F \colon (\m{C},P) \to (\m{D},Q) $ and $ K \colon (\m{C},P) \to (\m{E},R) $, there exists a unique morphism $ H \colon (\m{E},R) \to (\m{D},Q) $ such that $ H \circ K = F $. Part 2 says that under reasonable ``surjectivity'' assumptions on $ K $, there is a unique completion $ H $ iff $ \ker K \leq \ker F $. Our results on morphic images and the internal logic are used to give ``algebraic'' proofs of these results. Moreover, the fibered homomorphism theorem determines an orthogonal factorization system for $ \FA_{\m{L}} $.

	In the final section, we use the fibered homomorphism theorem to characterize two natural closure operators. The first operator characterizes the closure of a class of structures under the satisfaction of their common first-order theory. Taking submodels, homomorphic images and products does not preserve the satisfaction of first-order theories in the Tarskian semantics for classical first-order logic. However, for logics $ \m{L} $ given prop-categorical semantics in $ \FA_{\m{L}} $, the closure operator is characterized by $ \mathbb{H}\mathbb{S}\mathbb{P} $ where $ \mathbb{H} $, $ \mathbb{S} $ and $ \mathbb{P} $ close a collection of structures under natural notions of homomorphic image, submodel and product of structures in the prop-categorical semantics (Theorem \ref{thm: Characterization of Theory Closure.}). For example, if $ \m{L} $ is classical first-order logic, then the Tarskian $ Sg $-structures are essentially the morphisms $ F \colon (\m{C}_{Sg},P_{Sg}) \to (\Set,\s{P}) $, where $ \s{P} $ is the preimage functor, and $ (\m{C}_{Sg},P_{Sg}) $ is the classifying prop-category of the minimal $ Sg $-theory. Then $ \FA_{\m{L}}((\m{C}_{Sg},P_{Sg}),(\Set,\s{P})) $ provides complete semantics for $ \m{L} $ (restricted to the signature $ Sg $). Just as the two element Boolean algebra and the class of all Boolean algebras forms complete semantics for classical propositional logic, $ \FA_{\m{L}} $ also provides complete semantics for $ \m{L} $, where an $ Sg $-structure is a morphism $ F \colon (\m{C}_{Sg},P_{Sg}) \to (\m{C},P) $, for some $ (\m{C},P) \in \Ob(\FA_{\m{L}}) $. In this extended semantics, it is natural to look at the \textit{coslice category} $ (\m{C}_{Sg},P_{Sg}) \! \downarrow \! \FA_{\m{L}} $ of the underlying $ 1 $-category of $ \FA_{\m{L}} $ instead of $ \FA_{\m{L}}((\m{C}_{Sg},P_{Sg}),(\Set,\s{P})) $. It is in this category $ (\m{C}_{Sg},P_{Sg}) \! \downarrow \! \FA_{\m{L}} $ that we define the operations of $ \mathbb{H}, \mathbb{S} $ and $ \mathbb{P} $, and obtain the fibered HSP-result.

	The second operator closes a collection $ \m{X} \subseteq \Ob(\FA) $ under logical consequence: that is, the closure of $ \m{X} $ is the largest $ \m{Y} \supseteq \m{X} $ such that $ \vDash_{\m{Y}} \ = \ \vDash_{\m{X}} $, i.e.\ they define the same logic. Since $ \FA_{\m{L}} $ is stable under products and subprop-categories,  $ \m{Y} \supseteq \mathbb{S}\mathbb{P}(\m{X}) $, where here $ \mathbb{S} $ is closure under subprop-categories and $ \mathbb{P} $ is closure under products of prop-categories. It turns out that we need another operation which we call $ \mathbb{U} $, and Theorem \ref{thm: Characterization of Logic Closure.} shows that $ \m{Y} = \mathbb{U}\mathbb{S}\mathbb{P}(\m{X}) $. This result is a natural extension of a corresponding result by Blok and J\'onsson in \cite{Blok2006}, which shows that for a collection of algebras $ \m{V} $, closure under the equational consequence is given by $ \mathbb{U}_{\lambda} \mathbb{S}\mathbb{P}(\m{V}) $, where $ \lambda $ is the cardinality of the set of variables over which equations are defined, and $ B \in \mathbb{U}_{\lambda}(\m{V}) $ if every $ \lambda $-generated subalgebra $ A \leq B $ is in $ \m{V} $. (See Theorem \ref{thm: Fixed signature adjunction} for the directly analogous ``fibered'' result which considers a fixed first-order signature and uses an operator $ \mathbb{U}_{Sg} $.)

	\section{First-Order Logics} \label{sec: First-Order Logics}
	
	This section begins by presenting the first-order syntax using \textit{contexts}. Contexts allow us to unambiguously interpret well-formed logical expressions in prop-categories which are introduced in the next section. After presenting typed equational logic the section ends by defining theories, first-order logics and derivability.
	
	Throughout the paper, we fix a \textbf{first-order language} $ \s{L} $, consisting of a collection $ \s{L}_{q} $ of quantifier symbols and a collection $  \s{L}_{\omega} $ of propositional connectives each with a designated arity in $ \omega $. For each $n \in \omega$, let $ \s{L}_{n} $ denote the $ n$-ary propositional connectives in $ \s{L}_{\omega} $. We consider multi-sorted first-order logic, and so a \textbf{signature} $ Sg $ is a collection of sort symbols $ \sigma, \tau, \gamma, \ldots $, typed function symbols $ f \colon \sigma_{1}, \ldots, \sigma_{n} \to \tau $, (we write $ f:\tau $ for $ n = 0 $), and typed relation symbols $ R \subseteq \sigma_{1}, \ldots, \sigma_{n} $.\footnote{$f \colon \sigma_1, \ldots, \sigma_n \to \tau $ and $ R \subseteq \sigma_1, \ldots, \sigma_n $ are merely typing assignments and do not assert that $ f $ is a morphism or $R$ is a subset of some set.} (we write $ R \subseteq \emptyset $ for $ n = 0 $)

	To provide prop-categorical semantics, well-formed terms and formulas must be ``in context''. A \textbf{context} $ \Gamma = [x_{1}:\sigma_{1}, \ldots, x_{n}: \sigma_{n}] $ is a list of distinct variables, $ x_{1}, \ldots, x_{n} $ where each variable $ x_{i} $ is assigned a type (sort symbol) $ \sigma_{i} $. We let $ M : \tau \ [\Gamma] $ be the assertion that the term $ M $ is of type $ \tau $ and $ \Gamma $ is a valid context for $ M $. The (well-formed) terms-in-context are defined inductively via the following typing rules:
	$$
	\AxiomC{ \phantom{ $ M_1 : \tau \ [\Gamma]  $  } }
	\UnaryInfC{$x: \sigma \ [\Gamma, x:\sigma, \Gamma^{\prime}] $}
	\DisplayProof
	\hskip 1.5 em
	\AxiomC{$ M_{1}:\sigma_{1} \ [\Gamma] \quad \ldots \quad M_{n}:\sigma_{n} \ [\Gamma] $}
	\RightLabel{,}
	\UnaryInfC{$ f(M_{1}, \ldots, M_{n}):\tau \ [\Gamma]$}
	\DisplayProof
	$$ 
	for each variable $ x $ and function symbol $ f \colon \sigma_{1}, \ldots, \sigma_{n} \to \tau $ in $ Sg $. These are the only typing rules for terms, and from these rules alone, one may prove that the typing rule
	\begin{equation}\label{eqn: admissible substitution typing rule}
	\AxiomC{$M :\sigma \ [\Delta]$}
	\AxiomC{$ N :\tau \ [\Gamma, x:\sigma,\Gamma^{\prime}]  $}
	\BinaryInfC{$ N[M/x] : \tau \ [\Gamma,\Gamma^{\prime}] $}
	\DisplayProof
	\end{equation}
	is admissible, where $ N[M/x] $ denotes the operation of substituting each occurrence of $ x $ in $ N $ with $ M $ and the variable sort pairs in $ \Delta $ are contained in the set of those in $ \Gamma,\Gamma^{\prime} $. We denote this by $ \VS(\Delta) \subseteq \VS(\Gamma,\Gamma^{\prime}) $.
	
	An \textbf{equation-in-context} $ M_{1} = M_{2} : \tau \ [\Gamma] $ is well-formed if $ M_{i} : \tau \ [\Gamma] $ holds for $ i \in \{ 1,2 \} $. In the sequel, all the logics we consider are built over the usual (typed) equational derivation system. The rules of this system are as follows:
	\begin{center}\label{fig: equational logic}
		\textbf{Equational Logic.}
		$$
		\AxiomC{$ $}
		\UnaryInfC{$ M = M : \sigma \ [\Gamma] $}
		\DisplayProof
		\hskip 1.5em
		\AxiomC{$ M = M^{\prime} : \sigma \ [\Gamma]$}
		\UnaryInfC{$ M^{\prime} = M : \sigma \ [\Gamma] $}
		\DisplayProof
		\hskip 1.5em
		\AxiomC{$ M = M^{\prime} : \sigma \ [\Gamma] \quad M^{\prime} = M^{\prime \prime} : \sigma \ [\Gamma] $}
		\UnaryInfC{$ M = M^{\prime \prime} : \sigma \ [\Gamma] $}
		\DisplayProof
		$$ 
		$$
		\bottomAlignProof
		\AxiomC{$ M = M^{\prime} : \sigma \ [\Delta] \quad N = N^{\prime} : \tau \ [\Gamma, x : \sigma, \Gamma^{\prime}] $}
		\RightLabel{ $ \VS(\Delta) \subseteq \VS(\Gamma,\Gamma^{\prime}). $}
		\UnaryInfC{$ N[M/x] = N^{\prime}[M^{\prime}/x] : \tau \ [\Gamma,\Gamma^{\prime}]  $}
		\DisplayProof
		$$ 
	\end{center}
	
	For a formula $ \phi $, we let $ \phi : \propo \ [\Gamma] $ be the assertion that $ \phi \, [\Gamma] $ is a well-formed formula-in-context. The \textbf{atomic} well-formed formulas-in-context are constructed using the rules
	$$
	\AxiomC{$M_{1}: \sigma_{1} \ [\Gamma] \quad \ldots \quad M_{n}: \sigma_{n} \ [\Gamma]$} 
	\UnaryInfC{$ R(M_{1}, \ldots, M_{n}): \propo \ [\Gamma] $}
	\DisplayProof
	\hskip 1.5em
	\AxiomC{$M_{1}:\sigma \ [\Gamma] \quad M_{2}:\sigma \ [\Gamma]$}
	\RightLabel{,}
	\UnaryInfC{$M_{1} =_{\sigma} M_{2} : \propo \ [\Gamma]$}
	\DisplayProof
	$$
	for each $ n $-ary relation symbol $ R \subseteq \sigma_{1}, \ldots, \sigma_{n} $. Note that $ M_{1} =_{\sigma} M_{2} \ [\Gamma] $ is considered a formula-in-context whereas $ M_{1} = M_{2}:\sigma \ [\Gamma] $ is not. Thus we have two different notions of equality and for many logics, $ =_{\sigma} $ is coarser than $ = $. 
	
	The well-formed formulas-in-context are recursively constructed from the atomic ones via the rules
	$$
	\AxiomC{$\phi_{1} : \propo \ [\Gamma] \quad \ldots \quad \phi_{n} : \propo \ [\Gamma]$}
	\UnaryInfC{$\Diamond( \phi_{1}, \ldots, \phi_{n}): \propo \ [\Gamma]$}
	\DisplayProof
	\hskip 1.5em
	\AxiomC{$\phi : \propo  \ [\Gamma, x : \sigma] $}
	\RightLabel{,}
	\UnaryInfC{$ \Omega_{x:\sigma}(\phi) : \propo \ [\Gamma] $}
	\DisplayProof
	$$
	for each quantifier symbol $ \Omega \in \s{L}_{q} $, all $ n \in \omega $ and each $ n $-ary propositional connective $ \Diamond \in \s{L}_{n} $. We consider formulas-in-context equal up to $ \alpha $-equivalence, i.e.\ up to renaming of bound variables. (With the syntax above, we can define the bound variables in $ \phi \, [\Gamma] $, to be those which occur in $ \phi $ but not in $ \Gamma $, and the free variables to be all other variables in $ \phi $ and $ \Gamma $.) This allows us to define the action of simultaneous term substitution on ($ \alpha $-equivalence classes of) formulas-in-context that avoids variable capture which we denote by $ \phi[M_{1}/x_{1}, \ldots, M_{n}/x_{n}] \, [\Gamma^{\prime}] $. Defining capture avoiding substitution on $ \alpha $-equivalence classes of formulas is quite technical, but considering formulas and terms in context greatly simplifies the construction and allows one to do so without explicit mention of free and bound variables.
	
	The well-formed \textbf{sequents-in-context} are of the form 
    $$ \phi_{1}, \ldots, \phi_{n} \vdash \phi_{n+1} \ [\Gamma], $$ 
    where $ n \in \omega $ and $ \phi_{i}: \propo \ [\Gamma] $, for each $ i \leq n+1 $.\footnote{We model single-conclusion sequent calculi since intuitionistic-type logics are the most common categorical logics, though multi-conclusion calculi can also be modeled by adding an extra condition on the prop-categorical semantics, analogous to Condition 4 in the next section (intuitionistic-type logics include full Lambek calculus in
	particular and most logical systems can be expressed as extensions of it).} We call the sequents-in-context and equations-in-context \textbf{assertions} and let $ A_{Sg} $ denote the collection of all $ Sg $-assertions. A \textbf{theory} $ T $ is an ordered pair $ T = (\Sg(T),\A(T)) $, where $ \Sg(T) $ is a signature and $ \A(T)$ is a collection of $ \Sg(T) $-assertions. For each signature $ Sg $, we let $ \Th_{Sg} $ denote the complete lattice of $ Sg $-theories ordered by $ T_{1} \leq T_{2} $ if and only if $ \A(T_{1}) \subseteq \A(T_{2}) $.
	
	Now that we have defined the well-formed expressions, a \textbf{(first-order) logic}\footnote{In Section \ref{sec: The 2-Categorical View of the Prop-Categorical Semantics} we define an action by substitution which can be used to define a (first-order) logic as a \textit{structural} closure operator.} is a closure operator $ \m{L} $ on the lattice of all $ Sg $-theories for each signature $ Sg $.  In the sequel, we will also assume all logics we consider satisfy the rules of equational logic in Figure \ref{fig: equational logic}. If each Given $ T,T^{\prime} \in \Th_{Sg} $, we let $ T_{\m{L}} $ denote the $ \m{L} $-closure of $ T $ and $ T \vdash_{\m{L}} T^{\prime} $, denote the assertion that $ T^{\prime} $ is \textbf{derivable} from $ T $, i.e.\ that $ T^{\prime} \leq T_{\m{L}} $. We say a logic $ \m{L}_1 $ is \textbf{stronger} than a logic $ \m{L}_{2} $, (denoted $ \m{L}_1 \geq \m{L}_2 $) if for each theory $ T $, $ T_{\m{L}1} \geq T_{\m{L}_{2}} $.

	\section{Prop-Categorical Semantics}\label{sec: Prop-Categorical Semantics}
	
	In this section, we define the class of prop-categories that provide semantics for first-order logics and give examples. Structures in prop-categories and the notion of a structure satisfying a theory are then defined and it is shown how a class of prop-categories defines a first-order logic.
	
	A \textbf{prop-category} $ (\m{C},P) $, is a category $ \m{C} $ with designated finite products, and a contravariant functor $ P \colon \m{C} \to \Pos $, where $ \Pos $ is the category of all partially ordered sets and monotone maps. The objects of $ \m{C} $ interpret the sorts, the morphisms interpret terms-in-context and the posets $ P(c) $ for each $ c \in \Ob(\m{C}) $ interpret the formulas-in-context. We define $ \Ob(\FA) $ to be the collection of all prop-categories which additionally satisfy:
	\begin{enumerate}
		\item For each $ c \in \Ob(\m{C}) $, $ P(c) $ is an $ \s{L}_{\omega} $-algebra\footnote{By $ \s{L}_{\omega} $-algebra we mean any algebra in the signature $ \s{L}_{\omega} $. We do not assume these algebras satisfy any particular collection of equations, nor do we assume their operations are monotone.} and for all $ f \in \Mor(\m{C}) $, $ P(f) $ is an $ \s{L}_{\omega} $-algebra homomorphism. \label{con: Images of morphisms are homomorphis.}
		\item For each $ c \in \Ob(\m{C}) $, there is a designated element $ Eq_{c} \in P(c  \times c) $.
		\item For each $ c \in \Ob(\m{C}) $, and each $ \Omega \in \s{L}_{q} $, there is a natural transformation $ \Omega_{(\cdot),c} \colon UP(\blank \times c) \Rightarrow UP $, where $ U \colon \Pos \to \Set $ is the forgetful functor.  \label{con: Naturality of quantifier maps.}
		\item For each $ c \in \Ob(\m{C}) $, $ P(c) $ has a designated binary operation $ \con^{P(c)} $, and nullary operation $ e_{c} \in P(c) $ so that $ (P(c),\otimes^{P(c)},e_{c}) $ is a monoid.
		\item For all $  \Omega \in \s{L}_{q} $, and all $ b,c,d \in \Ob(\m{C}) $,
		\begin{equation}\label{eqn: Extra quantifier condition.}
		\Omega_{b,1} \circ P(\pi_{1}^{b,1}) = id_{P(b)} \quad \text{and} \quad \Omega_{b, c \times d} \circ P(a_{b,c,d}) = \Omega_{b,c} \circ \Omega_{b \times c,d}, 
		\end{equation}
		where $ a_{b,c,d} \colon b \times (c \times d) \to (b \times c) \times d $ is the change-in-product isomorphism.  \label{con: 5 of prop-cat}
		\item $ Eq_{1} = e_{1 \times 1} $, and for all $ c_{1},c_{2} \in \Ob(\m{C}) $, $ c = c_{1} \times c_{2} $,\label{con: 6 of prop-cat.}
		\begin{equation}\label{eqn: Equality condition.}
		Eq_{c} = P(\langle \pi_{1}^{c_{1},c_{2}}\pi_{1}^{c,c}, \pi_{1}^{c_{1},c_{2}}\pi_{2}^{c,c} \rangle )Eq_{c_{1}} \con P(\langle \pi_{2}^{c_{1},c_{2}}\pi_{1}^{c,c}, \pi_{2}^{c_{1},c_{2}}\pi_{2}^{c,c} \rangle)Eq_{c_{2}}.
		\end{equation}
	\end{enumerate}
	
	Conditions 1-3 are necessary for structures in $ (\m{C},P) $ to interpret equations and formulas in context. Condition 4 and the requirement that the codomain of $ P $ is $ \Pos $ are needed to interpret sequents\footnote{In AAL, a (logical) matrix $ \langle A,F \rangle $ consists of an algebra $ A $ in the propositional language and a subset $ F \subseteq A $ of ``truth values''\cite[p.\ 183]{Font2016}. One can develop a ``fibered'' matrix semantics by dropping Condition 4 on prop-categories and defining a fibered matrix as $ \langle (\m{C},P),\s{F} \rangle $ where $ (\m{C},P) $ is a prop-category and $ \s{F}  \subseteq \sqcup_{c \in \Ob(\m{C})}P(c) $.}. Conditions 5 and 6 are properties satisfied by the classifying prop-categories associated to certain theories and are not strictly necessary to define the prop-categorical semantics. However, they must be assumed to interpret structures as morphisms out of classifying prop-categories, a key ingredient in the proofs of our main results.
	
	\begin{rem}
		The interpretation of the universal and existential quantifiers as right and left adjoints of $ P(\pi_{1}^{a,b}) $ respectively, seen in Example \ref{exp: General prop-category example.}, is due to Lawvere \cite{Lawvere1969}. One can show that Condition 5 follows from these adjointness conditions. In \cite{Lawvere1970}, Lawvere also provides an adjoint interpretation of the equality predicate. For an account in the context of intuitionistic logic, see \cite[p.~190]{Jacobs1999}. We use Conditions 5 and 6 instead of the corresponding adjointness conditions, because they are more general and in a form that will be directly used to prove essential properties of the fibered semantics\footnote{Proving 6 from the adjointness conditions requires some work and is proved in \cite[p.\ 198]{Jacobs1999} for a slightly different setup.}. Moreover, in the spirit of AAL, we seek to provide semantics for as broad of a variety of logics as possible. However, for many important examples $ \s{L}_{q} = \{ \forall, \exists \} $ and the quantifiers and equality do satisfy these adjointness conditions. Thus, one may safely replace Conditions $ 5 $ and $ 6 $ with these adjointness conditions provided corresponding ``adjoint'' derivation rules are added to the minimal logic $ \s{L}^{m} $ in Section \ref{sec: Classifying Prop-categories and General Completeness Theorems} to preserve the completeness results therein.
	\end{rem}

	\begin{exmp}\label{exp: General prop-category example.}
		Let $ L $ be an algebra of signature $ \s{L}_{\omega} $ which is also a complete lattice and let $ \s{L}_{q} = \{ \forall, \exists \} $. Then the contravariant functor $ \Set(\blank,L) \colon \Set \to \Pos $ determines a prop-category in $ \Ob(\FA) $. The operations on $ L $ extend pointwise to logical operations on $ \Set(A,L) $ for all sets $ A $. For each set $ X $, let $ e_{X} = \top $, $ \otimes \coloneqq \wedge $ and
		$$
		Eq_{X}(x_{1},x_{2}) = \begin{cases}
		\top & x_{1} = x_{2}, \\
		\bot & x_{1} \neq x_{2}.
		\end{cases} 
		$$ 
		For sets $ X,Y $, the quantifiers $ \exists_{X,Y} $ and $ \forall_{X,Y} $ are defined as the left and right adjoints of $ \Set(\pi_{1}^{X,Y},L) $ respectively. That is,
		$$
		\exists_{X,Y}(R)(x) = \bigvee_{y \in Y} R(x,y), \quad \text{and} \quad  \forall_{X,Y}(R)(x) = \bigwedge_{y \in Y} R(x,y).
		$$
		In particular, $ \Set(\blank,2) $ is the prop-category whose semantics corresponds to the usual Tarskian semantics for classical first-order logic\footnote{Except here, sorts can be interpreted as the emptyset.}.
		
		Another significant example in this class is $ \Set(\blank, [0,1]) $, where for each set $ A $, the poset $ \Set(A,[0,1]) $ is known as the set of all \textbf{fuzzy sets on $ A $}. There is a vast literature on fuzzy set theory (for an introduction see \cite{Hajek1998}) which generalizes set theory by considering partial or probabilistic membership in sets.
		
		An important class of logical operations on fuzzy sets are $ T $-norms \cite[p.\ 4]{Klement2000} and dually $ S $-norms. A $ T $-norm is an operation $ \odot \colon [0,1]^{2} \to [0,1] $ which makes $ ([0,1], \odot,1) $ a partially ordered monoid. Canonical examples of $ T $-norms include the minimum, product and the \L ukasiewicz $ T $-norm which is defined for $ x,y \in [0,1] $ by $ x \odot y \coloneqq \max \{ 0,x+y-1 \} $. In \cite{Thiele2000}, to each  $ T $-norm $ \odot $, a quantifier $ \Omega^{\odot} $, is defined on sets $ A,B $ and $ R \colon A \times B \to [0,1] $, by
		$$
		\Omega^{\odot}_{A,B}(R) \coloneqq \inf_{U \in \fin(B)}\{ \bigodot_{u \in U} R(\blank, u)  \},
		$$
		where $ \fin(B) $ is the set of all finite subsets of $ B $. In particular, note that if $ \odot(a,b) = \min\{ a,b \} $ then $ \Omega^{\odot} = \forall $ and if $ \odot(a,b) = \max\{ a,b \} $ then  $ \Omega^{\odot} = \exists $. One may verify that for an arbitrary $ T $-norm $ \odot $, $ \Omega^{\odot} $ satisfies Conditions 3 and 5. Thus $ \Set(\blank,[0,1]) $ may be extended with additional quantifiers $ \Omega^{\odot} $ for each $ T $-norm $ \odot $. Moreover, $ \Omega^{\odot} $ may also be extended with $ S $-quantifiers corresponding to $ S $-norms as defined in \cite{Thiele2000}. If $ \odot $ is the usual product, then for $ R \colon 1 \times 2 \to [0,1] $, $ R(0,0) = R(0,1) = 1/2 $, $ \forall_{1,2}(R)(0) = \exists_{1,2}(R)(0) = 1/2 $, whereas $ \Omega^{\odot}_{1,2}(R)(0) = 1/4 $. Moreover, observe that in this example $ \odot $ may be used for $ \otimes $ instead of $ \min $.
	\end{exmp}
	
	The following is a toy example which suggests integral operators are a fruitful place to look for nonclassical quantifiers and metrics for nonstandard interpretations of the equality predicate:
	\begin{exmp}
		Let $ \m{K} $ be the category whose objects are the compact subsets of $  \mathbb{R}^{n} $, for some $ n \in \omega $, where each $ \mathbb{R}^n $ is equipped with the Euclidean metric and whose morphisms are the continuous functions. Then $ \m{C}(\blank, \mathbb{R})^{op} \colon \m{K}^{op} \to \Pos $ is a prop-category where $ \m{C}(A,\mathbb{R}) $ is the set of continuous real-valued functions on $ A $. Post-composition by any continuous map $ h \colon \mathbb{R}^{n} \to \mathbb{R} $ can interpret each $ n $-ary connective $ \diamond \in \s{L}_n $. For each $ A \in \Ob(\m{K}) $, let $ e_A $ be the constant zero function and for $ f,g \in \m{C}(A,\mathbb{R}) $, define $
		f \otimes g \coloneqq \sqrt{f^2 + g^2} $.
		Let $ Eq_{A} \colon A \times A \to \mathbb{R} $ the Euclidean metric and we can define a quantifier $ \int $, which for all objects $ A,B $ and $ p \in \m{C}(A \times B, \mathbb{R}) $, $ \int_{A,B}(p) \coloneqq \int_{B}p(x,y) \, \mathrm{d}y $.
		
		Then $ (A, \otimes^{P(A)},e_{A}) $ is a monoid. Moreover, Condition 6 is satisfied. If $ \mathbb{R}^0 $ is given the discrete measure, $ \int_{A,1} \circ P(\pi_{1}^{A,1}) = id_{P(A)} $ and from Fubini's Theorem, $ \int_{(\cdot),(\cdot)} $ also satisfies the second part of Condition 5. Thus $ (\m{C},P) \in \Ob(\FA) $. Note, in general, $ \int_{A,B} $ is not left or right adjoint to $ P(\pi_{1}^{A,B}) $. To see this, consider $ A = B = [0,1] $, $ r(x,y) \colon A \times B \to \mathbb{R} $, $ r(x,y) = y $. Then $ \int_{A,B}(r) = 1/2 $ whereas $ \exists_{A,B}(r) = 1 $ and $ \forall_{A,B}(r) = 0 $. Also, equality is not transitive. In fact, $ Eq_{A}(a_{1},a_{2}) $ becomes ``more true'' the farther $ a_{1} $ and $ a_{2} $ are apart and so is better thought of as asserting inequality.
	\end{exmp}

	Let $ Sg $ be a signature and $ (\m{C},P) \in \Ob(\FA) $. An \textbf{$ Sg $-structure $ S $ in $ (\m{C},P) $} is an assignment of an object $ \db{\sigma} $ to each sort symbol $ \sigma $, a morphism 
	$
	\db{f} \colon \db{\sigma_{1}} \times \cdots \times \db{\sigma_{n}} \to \db{\tau}
	$
	to each function symbol $ f \colon \sigma_{1}, \ldots, \sigma_{n} \to \tau $ and an element
	$
	\db{R} \in P(\db{\sigma_{1}} \times \cdots \times \db{\sigma_{n}})
	$
	to each relation symbol $ R \subseteq \sigma_{1}, \ldots, \sigma_{n} $ in $ Sg $. (To disambiguate between structures, we sometimes write $ S \db{\cdot} $ instead of $ \db{\cdot} $.) For a context $ \Gamma = x_{1}:\sigma_{1}, \ldots, x_{n}:\sigma_{n} $, we define $ \db{\Gamma} \coloneqq \db{\sigma_{1}} \times \cdots \times \db{\sigma_{n}} $, in particular if $ n= 0 $, $ \db{\Gamma} \coloneqq 1 $. Each term-in-context $ M : \tau \ [\Gamma] $ may be given a unique interpretation $ \db{M : \tau \ [\Gamma]} $ defined recursively as follows:
	\begin{align*}\label{eqn: inductive interpretation of terms}
	\db{x_{i} : \sigma_{i} \ [x_{1}: \sigma_{1}, \ldots, x_{n} : \sigma_{n}]} & \coloneqq \pi_{i} \colon \db{\sigma_{1}} \times \cdots \times \db{\sigma_{n}} \to \db{\sigma_{i}} \\
	\db{c : \tau \ [\Gamma] } & \coloneqq \db{c} \circ !_{\db{\Gamma}} \\
	\db{f(M_{1}, \ldots, M_{n}) : \tau \  [\Gamma] } & \coloneqq \db{f} \circ \langle \db{M_{1}[\Gamma]}, \ldots, \db{M_{n} [\Gamma] } \rangle,
	\end{align*}
	where $ \pi_{i} $ is the $ i $-th projection map and $ !_{\db{\Gamma}} $ denotes the unique morphism from $ \db{\Gamma} $ to $ 1 $. Note we often abbreviate $ \db{M : \tau \ [\Gamma]} $ as $ \db{M [\Gamma] } $ as the type $ \tau $ can be inferred from the term $ M $ and context $ \Gamma $. 
	
	$ S $ interprets formulas-in-context recursively as follows:
	\begin{align*}
	\db{R(M_{1}, \ldots, M_{n}) \, [\Gamma]} \coloneqq & P( \langle \db{M_{1}[\Gamma]}, \ldots, \db{M_{n}[\Gamma]} \rangle  )(\db{R}), \\
	\db{M_{1}=_{\tau} M_{2} \ [\Gamma]} \coloneqq & P(\langle \db{M_{1}[\Gamma]}, \db{M_{2}[\Gamma]}\rangle)(\db{Eq_{\db{\tau}}}), \\
	\db{\Diamond(\phi_{1}, \ldots, \phi_{n}) \, [\Gamma]} \coloneqq & \Diamond^{P(\db{\Gamma})}(\db{\phi_{1} \, [\Gamma]}, \ldots, \db{\phi_{n} \, [\Gamma]}), \\
	\db{\Omega_{x:\sigma}(\phi) \, [\Gamma] } \coloneqq & \Omega_{\db{\Gamma},\db{\sigma}} \circ P(a^{S}_{\Gamma,x:\sigma}) (\db{\phi \, [\Gamma,x:\sigma]}),
	\end{align*}
	where $ a^{S}_{\Gamma,x:\sigma} \colon \db{\Gamma} \times \db{x:\sigma} \to \db{\Gamma,x:\sigma}$ is the change-in-product isomorphism.\footnote{This handles the case where $ \Gamma = \emptyset $ and $ \db{\Gamma,x:\sigma} = \db{x:\sigma} \neq \db{\Gamma} \times \db{x:\sigma} = 1 \times \db{x:\sigma}. $}

	Then, $S $ \textbf{satisfies an equation-in-context $ M_{1} = M_{2} : \tau \ [\Gamma] $}, if 
	$$
	\db{M_{1}:\tau \ [\Gamma]} = \db{M_{2} : \tau \ [\Gamma]}
	$$
	and \textbf{$ S $ satisfies a sequent-in-context} $ \phi_{1}, \ldots, \phi_{n} \vdash \phi_{n+1} \ [\Gamma],$ if for $ n \geq 0 $,
	$$
	\db{\phi_{1}[\Gamma]}\con \ldots \con \db{\phi_{n}[\Gamma]} \leq \db{\phi_{n+1} [\Gamma]}.
	$$
	In particular, for $ n=0 $, if $e_{\db{\Gamma}} \leq \db{\phi_{n+1} [\Gamma]}.$
	We define the \textbf{theory of $ S $}, denoted $ \Th(S) $, to be the $ Sg $-theory where $ \A(\Th(S)) $ is the collection of all assertions satisfied by $ S $. We say \textbf{$ S $ satisfies $ T $} or \textbf{$ S $ is a $ T $-model} if $ T \leq \Th(S) $.
	
	\begin{exmp}
		In \cite{Mostowski1957}, Mostowski introduced a general class of quantifiers and studied classical first-order logic extended with these quantifiers. Instead of $ \Set(\blank,2) $, we consider the isomorphic powerset functor $ \s{P} \colon \Set \to \Pos $, which sends a set $ A $ to its powerset $ \s{P}(A) $ and a function $ f \colon A \to B $ to its preimage operator. A \textbf{Mostowski quantifier}\footnote{In Mostowski's original formulation is more restrictive and requires quantifiers to be cardinality invariant. Our presentation here is essentially that of \cite{Westerstahl2019}.} $ \Omega $, is specified for each set $ A $, by a subset $ \Omega(A) \subseteq \s{P}(A) $. Then for $ A,B $ sets and $ R \subseteq A \times B $, we define
		$$
		\Omega_{A,B}(R) \coloneqq \{  a \in A : \{ b \in B : (a,b) \in R \} \in \Omega(B) \}.
		$$
		Note that quantifiers $ \forall $ and $ \exists $ belong to the class, defined for a set $ A $, by $ \forall(A) \coloneqq \{ A \} $ and $ \exists(A) \coloneqq \s{P}(A) \setminus \{ \emptyset \} $ respectively. A nonstandard Mostowski quantifier of early interest \cite{Kaufmann1985} is the ``there exists uncountably many'' quantifier $ Q_{1} $, where for each set $ A $, $ Q_{1}(A) $ is the collection of all uncountable subsets of $ A $. 
		
		One may verify that $ (\Set, \s{P}) $ augmented with a Mostowski quantifier $ \Omega $ satisfies Condition (3) and so structures in $ (\Set,\s{P}) $ can interpret the language of classical logic augmented with $ \Omega $. However, in general, many such natural quantifiers, fail to satisfy Condition (5). The first part of Equation \ref{eqn: Extra quantifier condition.} requires that $ Q_{1}(1) = \{ 1 \} $, and the second part of Equation \ref{eqn: Extra quantifier condition.} requires for all sets $ A,B $ and $ R \subseteq A \times B $, that
		\begin{equation}\label{eqn: Compact condition 2 for Mostowski quantifiers.}
		R \in Q(A \times B) \iff \{ a \in A : R(a,y) \in Q(B) \} \in Q(A).
		\end{equation}
		One sees that the second part of Condition (5) also fails for $ Q_{1} $, by considering $ A = B = \mathbb{R} $ and $ R = \mathbb{R} \times \{ 0 \} $. 
		
		One can do away with Condition (5) by expanding the syntax to allow quantifiers to bind arbitrary contexts as follows
		$$
		\AxiomC{$\phi : \propo  \ [\Gamma, \Gamma^{\prime}] $}
		\RightLabel{.}
		\UnaryInfC{$ \Omega_{\Gamma^{\prime}}(\phi) : \propo \ [\Gamma] $}
		\DisplayProof
		$$
		Then a structure $ S $ interprets $ \Omega_{\Gamma^{\prime}}(\phi) \, [\Gamma] $ as
		$$
		\db{\Omega_{\Gamma^{\prime}}(\phi) \, [\Gamma]} \coloneqq \Omega_{\db{\Gamma},\db{\Gamma^{\prime}}}\circ P(a_{\Gamma,\Gamma^{\prime}}^{S})(\db{\phi \, [\Gamma,\Gamma^{\prime}]}).
		$$
		This extension of the syntax is conservative for existential and universal quantifiers, but not for the counting quantifiers. For example, if $ \Omega^{n} $ is interpreted as ``there exists exactly $ n $.'' Then for each $ \phi : \mathrm{prop} \ [x:\sigma,y:\tau] $,
		$$
		\Omega^{2}_{x:\sigma,y:\tau} \, \phi(x,y) \, [
		\, ] \equiv \Omega^{1}_{y:\tau} \Omega^{2}_{x:\sigma} \, \phi(x,y) \vee \Omega^{2}_{y:\tau} \Omega^{1}_{x:\sigma} \, \phi(x,y) \, [\, ].
		$$
		To recover the results in this paper, one must add the following isomorphism invariant condition on the quantifiers: For each isomorphism $ f \colon c \to d $, 
		$$ \Omega_{(\cdot),c} \cdot P(\blank \times f) = \Omega_{(\cdot), d}, $$
		which is a condition Mostowski originally imposed so that his quantifiers are bijective invariant.

		We also note that similar considerations hold for the equality predicate and Condition 6.
	\end{exmp}

	A subcollection $ \m{X} \subseteq \Ob(\FA) $, defines a logic $ \vDash_{\m{X}} $, where given $ Sg $-theories $ T $ and $ T^{\prime} $, $ T \vDash_{\m{X}} T^{\prime} $ if every $ T $-model $ S \in (\m{C},P) \in \ca{X} $ is also a $ T^{\prime} $-model. Given a logic $ \m{L} $, we let $ \Ob(\FA_{\m{L}}) $ be the collection of all prop-categories in $ \Ob(\FA) $ such that $ \m{L} $ is sound with respect to $ \vDash_{\Ob(\FA_{\m{L}})} $, that is, whenever $ T \vdash_{\m{L}} T^{\prime} $, then $ T \vDash_{\Ob(\FA_{\m{L}})} T^{\prime} $.

	\section{Classifying Prop-Categories and General Completeness Theorems}\label{sec: Classifying Prop-categories and General Completeness Theorems}
	
	In the subsequent sections, we suppose $ \s{L} $ has designated operation symbols $ e \in \s{L}_{0} $ and $ \otimes \in \s{L}_{2} $. Given a sufficiently rich theory $ T $, we will construct the \textit{classifying prop-category} $ (\m{C}_{T},P_{T}) \in \Ob(\FA) $, which is the first-order analogue of the Lindenbaum-Tarski (LT) algebras of propositional logic. As with LT-algebras, each classifying prop-category $ (\m{C}_{T},P_{T}) $ contains a \textit{generic} $ T $-model $ G \in (\m{C}_{T},P_{T}) $, i.e.\ $ \Th(G) = T $, which interprets each term and formula in context as its associated equivalence class obtained by ``quotienting out'' by the theory $ T $.
	
	We call a logic $ \m{L} $ over $ \s{L} $ \textbf{adequate}, if $ \m{L} $ is stronger than $ \m{L}^{m} $, where $ \m{L}^{m} $ is the first-order logic whose rules are listed in Figure \ref{fig: Rules for L^m.}.\footnote{In AAL, propositional logics which satisfy $ \Diamond $-Cong for each connective are called selfextensional, which is the weakest class of logics in the Fregean hierarchy. \cite{Wojcicki1979} \cite[p.\ 419]{Font2016} Kleene's strong 3-valued logic $ K3 $ \cite[p.\ 332]{Kleene1952} is a simple example of a non-selfextensional propositional logic. However, A first-order version of $ K3 $ can still be modeled as the fragment of an adequate logic by restricting to sequents with empty antecedent.} 
	
	\begin{figure}
		\begin{center}
			$$
			\AxiomC{}
			\RightLabel{Ax}
			\UnaryInfC{$\phi \vdash \phi \ [\Gamma]$}
			\DisplayProof
			\hskip 1.5 em
			\AxiomC{$ \phi \vdash \psi \ [\Gamma] \quad \psi \vdash \theta \ [\Gamma] $}
			\RightLabel{Cut}
			\UnaryInfC{$ \phi \vdash \theta \ [\Gamma]$}
			\DisplayProof
			$$
			
			$$
			\AxiomC{$ \Phi \vdash \psi \ [\Gamma] $}
			\RightLabel{Cwk}
			\UnaryInfC{$ \Phi \vdash \psi \ [\Gamma, x:\sigma] $}
			\DisplayProof
			\hskip 1.5em
			\AxiomC{$ M = M^{\prime} :\sigma \ [\Delta] \quad  \Phi \vdash \psi \ [\Gamma, x: \sigma, \Gamma^{\prime}] $}
			\RightLabel{Sub }
			\UnaryInfC{$ \Phi[M/x] \vdash \psi[M^{\prime}/x] \ [\Gamma, \Gamma^{\prime} ] $}
			\DisplayProof
			$$
			
			$$
			\AxiomC{$ \phi \dashv \vdash \psi \ [\Gamma,x:\sigma] $}
			\RightLabel{$\Omega$-Cong}
			\UnaryInfC{$ \Omega_{x:\sigma}(\phi) \vdash \Omega_{x:\sigma}(\psi) \ [\Gamma] $}
			\DisplayProof
			\hskip 1.5em
			\AxiomC{ $\phi_{1} \dashv \vdash \phi_{1}^{\prime} \ [\Gamma] \quad \ldots \quad \phi_{n} \dashv \vdash \phi_{n}^{\prime} \ [\Gamma] $}
			\RightLabel{$ \Diamond $-Cong}
			\UnaryInfC{$ \Diamond(\phi_{1}, \ldots, \phi_{n}) \vdash  \Diamond(\phi_{1}^{\prime}, \ldots, \phi_{n}^{\prime}) \ [\Gamma]$}
			\DisplayProof
			$$
			
			$$
			\AxiomC{$\Phi, \alpha, \beta,\Psi \vdash \theta \ [\Gamma] $}
			\RightLabel{$ \con $-Ref}
			\doubleLine
			\UnaryInfC{$ \Phi,\alpha \con \beta, \Psi \vdash \theta \ [\Gamma] $}
			\DisplayProof
			\hskip 1.5 em
			\AxiomC{$ \Phi, e, \Psi \vdash \phi \ [\Gamma] $}
			\RightLabel{$ e $-Ref}
			\doubleLine
			\UnaryInfC{$  \Phi,\Psi \vdash \phi \ [\Gamma]$}
			\DisplayProof
			$$
		\end{center}
		\caption{Rules for $ \m{L}^{m} $, where $ \VS(\Delta) \subseteq \VS(\Gamma,\Gamma^{\prime}) $ and $ \Omega \in \s{L}_{q} $, $ n \in \omega $ and $ \Diamond \in \s{L}_{n} $. Note: double lines mean the rule may be applied top-down or bottom-up.}\label{fig: Rules for L^m.}
	\end{figure}
	
	\begin{thrm}\label{thm: Existence of classifying category and generic model.} 
		If $ \m{L} $ is adequate, then for each $ \m{L} $-theory $ T $, one can construct a prop-category $ (\m{C}_{T},P_{T}) $, which contains a structure $ G $, such that $ T = \Th(G) $. 
	\end{thrm}
	\begin{proof}
	To construct $ \m{C}_{T} $, from an $ \m{L} $-theory $ T $, we consider contexts, formulas and terms in context equal up to renaming of their free variables. Morphisms are sequences of equivalence classes of terms in a common context:
	$$
	\gamma \colon \Gamma^{\prime} \to \Gamma = [N_{1}:\sigma_{1} \ [\Gamma^{\prime}] / \! \sim, \ldots, N_{n}: \sigma_{n} \ [\Gamma^{\prime}] / \! \sim ],
	$$
	where, $ \Gamma = x_{1}:\sigma_{1}, \ldots, x_{n}:\sigma_{n} $ and
	$$ N : \tau \ [\Gamma] \sim M : \tau \ [\Gamma] \iff N = M : \tau \ [\Gamma] \in T. $$
	Composition of morphisms is by component-wise substitution and one may show using the rules of equational logic that $ \m{C}_{T} $ is a category with strictly associative finite products \cite[p.\ 30]{Pitts2000}. 
	
	For each context $ \Gamma $, we say $ \phi \, [\Gamma] \sim \psi \, [\Gamma] $ if and only if $ \phi \vdash \psi \, [\Gamma], \psi \vdash \phi \, [\Gamma] \in T $. From $ \mathrm{Ax} $ and $ \mathrm{Cut} $, $ \sim $ defines and equivalence relation on formulas-in-context $ \Gamma $ and we let $ P_{T}(\Gamma) $ be the partial order of all such classes where  $ \phi \, [\Gamma]/ \! \sim \ \leq \, \psi \, [\Gamma]/ \! \sim $ if and only if $ \phi \vdash \psi \, [\Gamma] \in T$. For each morphism $ \gamma \colon \Gamma^{\prime} \to \Gamma $, we define $ P_{T}(\gamma) \colon P_{T}(\Gamma) \to P_{T}(\Gamma^{\prime}) $ by simultaneous substitution:
	$$
	P_{T}([N_{1}:\sigma_{1} \, [\Gamma^{\prime}] / \! \sim, \ldots, N_{n}: \sigma_{n} \, [\Gamma^{\prime}] / \! \sim ])(\phi \, [\Gamma]/\! \sim) \coloneqq \phi [N_{1}/x_{1}, \ldots, N_{n}/x_{n}] \, [\Gamma^{\prime}]/ \! \sim. 
	$$
	For each $ n \in \omega $ and $ \Diamond \in \s{L}_{n} $, we define:
	$$
	\Diamond^{\Gamma}(\phi_{1}[\Gamma]/ \! \sim, \ldots, \phi_{n}[\Gamma]/ \! \sim) \coloneqq \Diamond(\phi_{1}, \ldots, \phi_{n})[\Gamma]/ \! \sim.
	$$
	For each context $ \Gamma $, define $ e_{\Gamma} \coloneqq e \, [\Gamma] / \! \sim $. For each quantifier symbol $ \Omega \in \s{L}_{q} $, and $ \phi \, [\Gamma,\Gamma^{\prime}] $, where $ \Gamma^{\prime} = y_{1}:\tau_{1}, \ldots, y_{m}:\tau_{m} $, we interpret $ \Omega_{\Gamma,\Gamma^{\prime}} $ as
	$$
	\Omega_{\Gamma,\Gamma^{\prime}}(\phi \, [\Gamma,\Gamma^{\prime}]/ \! \sim) \coloneqq \Omega_{y_{1}:\tau_{1}} \ldots \Omega_{y_{m}:\tau_{m}}(\phi)\, [\Gamma] / \! \sim.
	$$
	For each $ \Gamma = x_{1}:\sigma_{1}, \ldots, x_{n}:\sigma_{n} $, let $ \Gamma^{\prime} = x_{1}^{\prime}:\sigma_{1}, \ldots, x_{n}^{\prime}, \sigma_{n} $. We define $$
	Eq_{\Gamma} \coloneqq x_{1}=_{\sigma_{1}}x_{1}^{\prime} \con \ldots \con x_{n} =_{\sigma_{n}} x_{n}^{\prime} \ [\Gamma,\Gamma^{\prime}]/ \! \sim.
	$$
	In particular, when $ \Gamma = [ \, ]$, $
	Eq_{[ \, ]} = e [ \, ]/ \! \sim \, = e_{[ \, ] \times [ \, ]}.$
	Using the rules in Figure \ref{fig: Rules for L^m.} one shows that these operations are well-defined and $ (\m{C}_{T},P_{T}) \in \Ob(\FA) $.
	
	The generic $ T $-model $ G  \in (\m{C}_{T},P_{T}) $, is defined on sorts by $ G \db{\sigma} \coloneqq x:\sigma $, on function symbols $ f \colon \sigma_{1}, \ldots, \sigma_{n} \to \tau $, by $ G\db{f} \coloneqq [f(\bar{x}) : \tau \ [\Gamma] / \! \sim] $ and on relation symbols $ R \subseteq \sigma_{1}, \ldots, \sigma_{n} $ by $ G \db{R} \coloneqq R(\bar{x}) \, [\Gamma] / \! \sim $. By induction one shows that $ G\db{M : \tau \ [\Gamma]} = [M:\tau \ [\Gamma]/ \! \sim] $ and $ G\db{\phi \, [\Gamma]} = \phi \, [\Gamma]/ \! \sim $. From $ \otimes $-Ref and $ e $-Ref it follows that $ \Th(G) = T  $.
\end{proof}
	We call the prop-category $ (\m{C}_{T},P_{T})$ constructed in the proof of Theorem \ref{thm: Existence of classifying category and generic model.} the \textbf{classifying prop-category of $ T $} and the $ T $-model $ G $ in the proof the \textbf{generic $ T $-model}.
	
	\begin{thrm}\label{thm: FAsL is a complete semnatics for Lm.}
		The logic $ \vDash_{\m{K}} $ provides a complete semantics for $ \vdash_{\m{L}^{m}} $, where $ \m{K} $ is any subclass of $ \Ob(\FA) $ containing $ \{ (\m{C}_{T},P_{T}): \text{$ T $ is an $ \m{L}^{m} $-theory} \} $.
	\end{thrm}
	\begin{proof}
		To show that $ \vdash_{\m{L}^{m}} \; \leq \ \vDash_{\Ob(\FA)} $, one verifies that each $ S \in (\m{C},P) \in \Ob(\FA) $ satisfies every interpretation of the rules defining $ \m{L}^{m} $. Now suppose $ T \vDash_{\m{K}} T^{\prime} $. Then the generic $ T_{\m{L}^{m}} $-model satisfies $ T^{\prime} $, and so $ T \vdash_{\m{L}^{m}} T^{\prime} $.
	\end{proof}
	As a consequence, we let $ (\m{C}_{Sg},P_{Sg}) \coloneqq (\m{C}_{Sg_{\m{L}^{m}}},P_{Sg_{\m{L}^{m}}})$. In $ (\m{C}_{Sg},P_{Sg}) $, morphisms $ \gamma \colon \Gamma \to \Gamma^{\prime} $ are just lists of terms and the elements of $ P_{Sg}(\Gamma) $ are just formulas, both up to $ \alpha $-equivalence and change in the name of free variables, which are equivalences considered in mathematical practice. Thus, $ (\m{C}_{Sg},P_{Sg}) $ is a suitable ``fibered'' analogue of the formula algebra from propositional algebraic logic.
	
	Suppose $ \m{L} $ is adequate and let $ \m{K} = \{ (\m{C}_{T},P_{T}): T \text{ is an $ \m{L} $-theory} \} $. Then using generic models, one proves  $ \vDash_{\m{K}}  \ \leq \ \vdash_{\m{L}} $. All that remains to prove completeness is to show soundness, i.e.\ $ \vdash_{\m{L}} \ \leq \ \vDash_{\m{K}} $. If $ \m{L} $ is defined by a sequent calculus, we only need to show each structure in a classifying prop-category satisfies the defining rules of $ \m{L} $. We have the following general result:
	\begin{thrm}\label{thm: General Completeness Theorem}
		Suppose $ \s{L}_{q} = \{ \forall, \exists \} $ and $ \m{L} $ is an extension of $ \m{L}^{m} $ possibly by  $ = $-Adj, $ \exists $-$ \mathrm{Adj} $, $ \forall $-$ \mathrm{Adj} $ shown in Figure \ref{fig: Adjoint rules.} and any number of structural and propositional connective rules\footnote{Propositional connective rules, such as modus ponens, are sequent rules whose meta-formulas do not include equality or any quantifiers. Structural rules, such as $ \mathrm{Cut} $, additionally do not include propositional connectives.}. Then $ \vdash_{\m{L}} \ = \ \vDash_{\m{K}} $, where $ \m{K} = \{ (\m{C}_{T},P_{T}): T \text{ is an $ \m{L} $-theory} \} $.
	\end{thrm}
	\begin{proof}
		Showing that $ \vDash_{\m{K}} $ satisfies the structural and propositional connective rules of $ \m{L} $ is straightforward.
		
		Suppose $ \m{L} $ satisfies $ \exists $-Adj and let
		$$
		\phi_{1}, \ldots, \phi_{n}, \psi \vdash \theta \ [\Gamma,x:\sigma],
		$$
		be an interpretation of the top line of $ \exists $-Adj. Let $ S \in (\m{C}_{T},P_{T}) \in \m{K} $ and for each $ i \in \{ 1, \ldots,n \} $, let $ S\db{\phi_{i} \, [\Gamma]} = \phi_{i}^{S} \, [\Gamma_{S}] / \! \sim $, $ S \db{\theta \, [\Gamma]} = \theta^{S} \, [\Gamma_{S}] \! \sim $ and $ S \db{\psi \, [\Gamma,x:\sigma]} = \psi^{S} \, [\Gamma_{S},\Gamma_{\sigma}]/ \! \sim $, where $ S \db{\sigma} = \Gamma_{\sigma} = x_{1}:\sigma_{1}, \ldots, x_{n}:\sigma_{n} $. Then,
		\begin{align*}
		& \bigcon_{i=1}^{n} S\db{\phi_{i}[\Gamma, x:\sigma] } \con S \db{\psi [\Gamma, x:\sigma]} \leq S \db{\theta [\Gamma, x: \sigma]} \\
		\iff \ & \bigcon_{i=1}^{n} G \db{\phi_{i}^{S}[\Gamma_{S},\Gamma_{\sigma}]} \con G \db{\psi^{S} [\Gamma_{S}, \Gamma_{\sigma}]} \leq G \db{\theta^{S}[\Gamma_{S},\Gamma_{\sigma}]}  \\
		\iff \ & T \vdash_{\m{L}} \phi_{1}^{S}, \ldots, \phi_{n}^{S}, \psi^{S} \vdash \theta^{S} \ [\Gamma_{S},\Gamma_{\sigma}] \\
		\iff \ & T \vdash_{\m{L}} \phi_{1}^{S}, \ldots, \phi_{n}^{S}, \exists_{x_{1}:\sigma_{1}} \ldots \exists_{x_{n}:\sigma_{n}}(\psi^{S}) \vdash \theta^{S} \ [\Gamma_{S}] \\
		\iff \ & \bigcon_{i=1}^{n} S \db{\phi_{i}[\Gamma]} \con S \db{\exists_{x:\sigma}(\psi)[\Gamma]} \vdash S \db{\theta [\Gamma]}.
		\end{align*}
		It follows that $ \vDash_{\m{K}} $ satisfies $ \exists $-Adj. The proof that $ \vDash_{\m{K}} $ satisfies $ \forall $-Adj is similar.

		Now suppose $ \m{L} $ satisfies $=$-Adj and let $ S $ be a structure in $ (\m{C}_{T},P_{T}) $. Consider an interpretation of $ = $-Adj:
		$$
		\AxiomC{$ \phi_{1}, \ldots, \phi_{n} \vdash \psi[x/x^{\prime}] \ [\Gamma,x:\sigma]$}
		\RightLabel{.}
		\doubleLine
		\UnaryInfC{$ \phi_{1}, \ldots, \phi_{n}, x =_{\sigma} x^{\prime} \vdash \psi \ [\Gamma, x:\sigma,x^{\prime}:\sigma]$}
		\DisplayProof
		$$
		Let $ S \db{\phi_{i}[\Gamma, x:\sigma]} = \phi_{i}^{S}[\Gamma_{S},\Gamma_{\sigma}] / \! \sim $, where $ \Gamma_{S} \in S \db{\Gamma} $, $ \Gamma_{\sigma} \in S\db{\sigma} $ and $ S\db{\psi [\Gamma, x:\sigma, x^{\prime},\sigma^{\prime}]} = \psi^{S} [\Gamma_{S},\Gamma_{\sigma},\Gamma_{\sigma}^{\prime}]/ \! \sim $, where $ \Gamma_{\sigma}^{\prime} \in S \db{\sigma} $. From $ \mathrm{Sub} $ and $ \mathrm{Cwk} $, one may derive the rule
		\begin{equation}\label{eqn: Context permutation}
		\AxiomC{$\Phi \vdash \psi \ [\Gamma, x:\sigma, y:\tau, \Gamma^{\prime}]$}
		\UnaryInfC{$\Phi \vdash \psi \ [\Gamma, y:\tau, x:\sigma, \Gamma^{\prime}]$}
		\DisplayProof
		\end{equation}
		Let $ \Gamma_{\sigma} = x_{1}:\sigma_{1}, \ldots, x_{n}:\sigma_{n} $ and $ \Gamma_{\sigma}^{\prime} = x_{1}^{\prime}:\sigma_{1}, \ldots, x_{n}^{\prime}:\sigma_{n} $. Let $ \bar{x} = x_{1},\ldots, x_{n} $ and $\bar{x}^{\prime} = x_{1}^{\prime}, \ldots, x_{n}^{\prime} $ be the lists of variables in $ \Gamma_{\sigma} $ and $ \Gamma_{\sigma}^{\prime} $ respectively. Then $ S \db{\psi[x/x^{\prime}][\Gamma,x:\sigma]} = \psi^{S}[\bar{x}/\bar{x}^{\prime}][\Gamma_{S},\Gamma_{\sigma}]/ \! \sim $ and
		\begin{align*}
		& S \db{\phi_{1}[\Gamma,x:\sigma]}\con \ldots \con S\db{\phi_{n}[\Gamma,x:\sigma]} \leq S \db{\psi[x/x^{\prime}][\Gamma,x:\sigma]} \\
		\iff \ & G\db{\phi_{1}^{S}[\Gamma_{S},\Gamma_{\sigma}]} \con \ldots \con G\db{\phi_{n}^{S}[\Gamma_{S},\Gamma_{\sigma}]} \leq G \db{\psi^{S}[\bar{x}/\bar{x}^{\prime}][\Gamma_{S},\Gamma_{\sigma}]} \\
		\iff \ & T \vdash_{\m{L}} \phi_{1}^{S}, \ldots, \phi_{n}^{S} \vdash \psi^{S}[\bar{x}/\bar{x}^{\prime}] \ [\Gamma_{S},\Gamma_{\sigma}]  \\
		\iff  \ & T \vdash_{\m{L}} \phi_{1}^{S}, \ldots, \phi_{n}^{S}, x_{1}=_{\sigma_{1}} x_{1}^{\prime}, \ldots, x_{n} =_{\sigma_{n}} x_{n}^{\prime} \vdash \psi^{S} \ [\Gamma_{S},\Gamma_{\sigma},\Gamma_{\sigma}^{\prime}] \tag{ $=$-Adj \text{and} \ref{eqn: Context permutation}} \\
		\iff \ & T \vdash_{\m{L}} \phi_{1}^{S}, \ldots, \phi_{n}^{S}, \bigcon_{i=1}^{n}x_{i}=_{\sigma_{i}}x_{i}^{\prime} \vdash \psi^{S} \ [\Gamma_{S}, \Gamma_{\sigma},\Gamma_{\sigma}^{\prime}] \tag{$\otimes$-Ref or $ e $-Ref}   \\
		\iff \ & S \db{\phi_{1}[\Gamma, x:\sigma,x^{\prime}:\sigma]} \con \ldots \con S\db{\phi_{n}[\Gamma,x:\sigma, x^{\prime}:\sigma]} \con S \db{x=_{\sigma}x^{\prime}[\Gamma,x:\sigma,x^{\prime}:\sigma]} \\ 
		& \leq S \db{\psi [\Gamma,x:\sigma,x^{\prime}:\sigma]}.
		\end{align*}
		Therefore, $ S $ satisfies all interpretations of $ = $-Adj.
	\end{proof}
	
	\begin{figure}
		\begin{center}
			$$
			\AxiomC{$ \phi_{1}, \ldots, \phi_{n} \vdash \psi[x/x^{\prime}] \ [\Gamma,x:\sigma]$}
			\RightLabel{$ = $-Adj}
			\doubleLine
			\UnaryInfC{$ \phi_{1}, \ldots, \phi_{n}, x =_{\sigma} x^{\prime} \vdash \psi \ [\Gamma, x:\sigma,x^{\prime}:\sigma]$}
			\DisplayProof
			$$
			
			$$
			\AxiomC{$\Phi \vdash \psi \ [\Gamma, x:\sigma]$}
			\doubleLine
			\RightLabel{$ \forall$-$ \mathrm{Adj} $}
			\UnaryInfC{$ \Phi,  \vdash \forall_{x,\sigma}(\psi) \ [\Gamma] $}
			\DisplayProof
			\hskip 1.5 em
			\AxiomC{$\Phi, \psi \vdash \theta \ [\Gamma, x:\sigma]$}
			\doubleLine
			\RightLabel{$ \exists$-$ \mathrm{Adj} $}
			\UnaryInfC{$ \Phi, \exists_{x:\sigma}(\psi) \vdash \theta \ [\Gamma] $}
			\DisplayProof
			$$
		\end{center}
		\caption{Adjoint rules for $ \forall $, $ \exists $ and $ = $.}\label{fig: Adjoint rules.}
	\end{figure}

	\section{The 2-Categorical View of the Prop-Categorical Semantics}\label{sec: The 2-Categorical View of the Prop-Categorical Semantics}
	
	We extend $ \Ob(\FA) $ to a $ 2 $-category by defining morphisms $ F \colon (\m{C},P) \to (\m{D},Q) $ by the following data:
	(1) a product preserving functor $ F^{o} \colon \mathcal{C} \to \mathcal{D} $ and (2) a natural transformation $ F^{p} \colon P \Rightarrow Q \circ F^{o} $,
	\begin{equation*}
	\adjustbox{scale=1.1,center}{
		\begin{tikzcd}
		\m{C} \arrow[rr, "P"{name = {P}}] \ar[swap]{dr}{F^{o}} & & \Pos	\\
		& \m{D}  \ar[swap]{ru}{Q}  \arrow[Rightarrow,from = P, shorten = 5, swap, "F^{p}"] & 
		\end{tikzcd}}
	\end{equation*}
	such that the following hold:
	\begin{enumerate}
		\item For each $ c \in \Ob(\m{C}) $, $ F^{p}_{c} \colon P(c) \to Q\circ F^{o}(c) $ is an $ \s{L}_{\omega} $-algebra homomorphism.
		\item For all $ \Omega \in \s{L}_{q} $ and all $ b,c \in \Ob(\m{C}) $, 
		\begin{equation*}
		F^{p}_{b} \circ \Omega_{b,c} = \Omega_{F^{o}b,F^{o}c} \circ Q(a_{F,b,c}^{-1}) \circ F^{p}_{b \times c}.
		\end{equation*}
		\item For all $ c \in \Ob(\m{C}) $,
		$$ F^{p}_{c \times c}(Eq_{c}) = Q(a_{F,c,c})(Eq_{F^{o}c}), $$
	\end{enumerate}
	where $ a_{F,b,c} \colon F^{o}(b \times c) \to F^{o}b \times F^{o}c $ is the change-in-product isomorphism. Given $ K \colon (\m{D},Q) \to (\m{E},R) $, we define $ K \circ F $ by $ (K \circ F)^{o} \coloneqq K^{o} \circ F^{o} $ and $ (K \circ F)^{p} \coloneqq K^{p}_{F^{o}} \cdot F^{p}  $. For parallel morphisms $ F,H \colon (\m{C},P) \to (\m{D},Q) $,  we define a $ 2 $-cell $ \eta \colon F \Rightarrow H $ to be a natural transformation $ \eta \colon F^{o} \Rightarrow H^{o} $ such that $ F^{p} = Q\eta \cdot H^{p} $.
	\begin{equation*}
	\adjustbox{scale=1.1,center}{
		\begin{tikzcd}
		\m{C} \arrow[rr, "P"{name = {P}}]  \arrow[dr, "H^{o}"{name = H^{o}}, bend right = 65,swap] \arrow[dr, "F^{o}"{name = F^{o}}] \arrow[Rightarrow, from = F^{o}, to = H^{o}, shorten = 5, swap, "\eta"] & & \Pos & P \arrow[Rightarrow]{d}{F^{p}} \arrow[Rightarrow]{dr}{H^{p}} &	\\
		& \m{D}  \arrow[swap]{ru}{Q} \arrow[Rightarrow,from = P, bend left = 40,pos = .53, shorten = 5, "H^{p}"] \arrow[Rightarrow, bend right = 15,from = P, shorten = 5, "F^{p}"]  &  & QF^{o} & QH^{o} \arrow[Rightarrow]{l}{Q\eta}
		\end{tikzcd}}
	\end{equation*}
	Compositions of $ 2 $-cells is just as it is in $ \Cat $. If $ \eta $ is a natural isomorphism, we call $ \eta $ a $ 2 $-isomorphism. For a logic $ \m{L} $, we now let $ \FA_{\m{L}} $ be the corresponding full sub-2-category of $ \FA $.

	\begin{thrm}\label{thm: 1 morphisms preserve models.}
		Let $ (\m{C},P),(\m{D},Q) \in \Ob(\FA) $ and $ F \colon (\m{C},P) \to (\m{D},Q) $ be a morphism. Then for each theory $ T $, each $ T $-model $ S $ in $ (\m{C},P) $ gives rise to a $ T $-model $ F(S) $ in $ (\m{D},Q) $.
	\end{thrm}
	
	\begin{proof}
		First we define the structure $ F(S) $. For each sort symbol $ \sigma $, $ F(S)\db{\sigma} \coloneqq F^{o}(S\db{\sigma}) $, for each function symbol $ f \colon \sigma_{1}, \ldots, \sigma_{n} \to \tau $, $ F(S)\db{f} \coloneqq F^{o}(S\db{f}) \circ a_{\Gamma}^{-1} $, and for each relation symbol $ R \subseteq \sigma_{1}, \ldots, \sigma_{n} $, $ F(S)\db{R} \coloneqq Q(a_{\Gamma}^{-1}) \circ F^{p}_{S\db{\Gamma}}(S\db{R}) $, where $ a_{\Gamma} \colon F^{o}(S \db{\Gamma} ) \to F(S)\db{\Gamma} $ is the change-of-product morphism and $ \Gamma = x_{1}:\sigma_{1}, \ldots, x_{n}:\sigma_{n} $. A product preserving functor preserves the satisfaction of equations-in-context since by induction on the complexity of an arbitrary term-in-context $ M : \tau \ [\Gamma] $,
		$$
		F(S)\db{M : \tau \ [\Gamma]} = F^{o}(S \db{M : \tau \ [\Gamma]}) \circ a_{\Gamma}^{-1}.
		$$
		If, for each formula-in-context $ \phi \, [\Gamma] $,
		\begin{equation}\label{eqn: goal to prove morphisms preserve modesl}
		F(S)\db{\phi [\Gamma]} = Q(a^{-1}_{\Gamma}) (F^{p}_{S\db{\Gamma}  } (S \db{\phi [\Gamma]}  )),
		\end{equation}
		then if $ S $ satisfies a sequent-in-context $ \phi_{1}, \ldots, \phi_{n} \vdash \psi \ [\Gamma] $, since $ F^{p}_{S\db{\Gamma}} $ and $ Q(a^{-1}) $ are monotone $ \s{L} $-algebra homomorphisms,
		\begin{align*}
		& \bigcon_{i=1}^{n}S\db{\phi_{i} [\Gamma]} \leq S \db{\psi [\Gamma]} \\
		\implies \ & \bigcon_{i=1}^{n} Q(a_{\Gamma}^{-1}) \circ F^{p}_{S \db{\Gamma}} (S \db{\phi_{i} [\Gamma]} ) \leq Q(a_{\Gamma}^{-1}) \circ F^{p}_{S \db{\Gamma} } (S\db{\psi [\Gamma]}) \\
		\implies \ & \bigcon_{i = 1}^{n} F(S)\db{\phi_{i} [\Gamma]} \leq F(S) \db{\psi [\Gamma]},
		\end{align*}
		where if $ n= 0 $, then $ \bigcon_{i=1}^{n} = e_{S\db{\Gamma}} $. And so $ F(S) $ satisfies $ \phi_{1}, \ldots, \phi_{n} \vdash \psi \ [\Gamma] $. It follows that $ F(S) $ is a $ T $-model and so we only need to prove that Equation \ref{eqn: goal to prove morphisms preserve modesl} holds.
		
		Consider a formula-in-context $ \phi \, [\Gamma] $. Suppose $ \phi \, [\Gamma] $ is atomic. Case 1: $ \phi \, [\Gamma] $ is of the form $ R(M_{1}, \ldots, M_{n})[\Gamma] $, where $ R \subseteq \tau_{1}, \ldots, \tau_{n} $, and each $ M_{i}:\tau_{i} \,  [\Gamma] $ is well-formed. If $ \gamma $ is the list $ \gamma = [M_{1}[\Gamma], \ldots, M_{n}[\Gamma]]$, we define $$
		S\db{\gamma} \coloneqq \langle S\db{M_{1}[\Gamma]}, \ldots, S\db{M_{n}[\Gamma]} \rangle.
		$$ 
		Note that 
		$$
		F(S)\db{\gamma} = a_{\Gamma^{\prime}} \circ F^{o}(S \db{\gamma})\circ a_{\Gamma}^{-1}
		$$
		Then,
		\begin{align*}
		& F(S) \db{R(M_{1}, \ldots,M_{n})[\Gamma]} \\
		= \ & Q(F(S)\db{\gamma})(F(S)\db{R}) \\
		= \ & Q(F(S)\db{\gamma})\circ Q(a_{\Gamma^{\prime}}^{-1})(F^{p}_{S\db{\Gamma^{\prime}}}(S\db{R})) \\
		= \ & Q(a_{\Gamma}^{-1}) \circ Q(F^{o}(S\db{\gamma})) \circ F^{p}_{S \db{\Gamma^{\prime}}}(S \db{R}) \\
		= \ & Q(a_{\Gamma}^{-1}) \circ F^{p}_{S\db{\Gamma}} \circ P(S\db{\gamma})(S\db{R}) \\
		= \ & Q(a_{\Gamma}^{-1})\circ F^{p}_{S \db{\Gamma}}(S \db{R(M_{1}, \ldots,M_{n})[\Gamma]}).
		\end{align*}
		Otherwise, $ \phi \, [\Gamma] $ is of the form $ M_{1}=_{\tau} M_{2} \ [\Gamma] $. Let $ \Gamma^{\prime} = [x_{1}:\tau,x_{2}:\tau] $ and $ \gamma \colon \Gamma \to \Gamma^{\prime} $ be the context morphism represented by $ [M_{1}[\Gamma],M_{2}[\Gamma]] $. Then,
		\begin{align*}
		& F(S) \db{M_{1}=_{\tau}M_{2}[\Gamma]} \\
		= \ & Q(F(S)\db{\gamma})(Eq_{F(S)\db{\tau}}) \\
		= \ & Q(a_{\Gamma}^{-1}) \circ Q(F^{o}(S\db{\gamma})) \circ Q(a_{\Gamma^{\prime}}) \circ Eq_{F(S)\db{\tau}} \\
		= \ & Q(a_{\Gamma}^{-1}) \circ Q(F^{o}(S \db{\gamma})) \circ F^{p}_{S \db{\tau} \times S\db{\tau}}(Eq_{S\db{\tau}}) \\
		= \ & Q(a_{\Gamma}^{-1}) \circ F^{p}_{S \db{\Gamma}} \circ P(S \db{\gamma})(Eq_{S\db{\tau}}) \\
		= \ & Q(a_{\Gamma}^{-1}) \circ F^{p}_{S \db{\Gamma}}(S \db{M_{1}=_{\tau}M_{2}[\Gamma]}).
		\end{align*}
		Now suppose $ \phi \, [\Gamma] $ is of the form $ \Diamond(\phi_{1}, \ldots, \phi_{n}) [\Gamma] $, where $ \Diamond $ is some $ n $-ary operation in $ \s{L}_{\omega} $ and for each $ i \in \{ 1, \ldots,n \} $, $ F(S)\db{\phi_{i}[\Gamma]} = Q(a^{-1}_{\Gamma}) \circ F^{p}_{S \db{\Gamma}}(S\db{\phi_{i}[\Gamma]}) $. Then,
		\begin{align*}
		& F(S) \db{\Diamond(\phi_{1}, \ldots, \phi_{n})[\Gamma]} \\
		= \ & \Diamond^{Q(F(S)\db{\Gamma})}(F(S)\db{\phi_{1}[\Gamma]}, \ldots, F(S)\db{\phi_{n}[\Gamma]}) \\
		= \ & \Diamond^{Q(F(S)\db{\Gamma})}(Q(a_{\Gamma}^{-1})\circ F^{p}_{S \db{\Gamma}}(S \db{\phi_{1}[\Gamma]}), \ldots,Q(a_{\Gamma}^{-1})\circ F^{p}_{S \db{\Gamma}}(S \db{\phi_{n}[\Gamma]}) ) \\
		= \ & Q(a_{\Gamma}^{-1}) \circ F^{p}_{S \db{\Gamma}}(\Diamond^{P(S \db{\Gamma})}(S\db{\phi_{1}[\Gamma]}, \ldots, S \db{\phi_{n}[\Gamma]}   )) \\
		= \ & Q(a_{\Gamma}^{-1}) \circ F^{p}_{S \db{\Gamma}}(S \db{\Diamond(\phi_{1}, \ldots,\phi_{n})[\Gamma]}).
		\end{align*}
		Suppose $ \phi \, [\Gamma] $ is of the form $ \Omega_{x:\sigma}(\psi) \, [\Gamma] $ for $ \Omega \in \s{L}_{q} $. Let $ I = S \db{\Gamma} $ and $ X = S\db{\sigma} $. For $ c = c_{1} \times \cdots \times c_{n} $, we let $ a_{c} \colon F^{o}(c) \to F^{o}c_{1} \times \cdots \times F^{o}c_{n} $, and $ a_{c_{1},c_{2}} \colon F^{o}(c_{1} \times c_{2}) \to F^{o}c_{1} \times F^{o}c_{2} $ be the change in product isomorphisms.
		Then \begin{align*}
		& F(S) \db{\Omega_{x:\sigma}(\psi)[\Gamma]} \\
		= \ & \Omega_{F(S)\db{\Gamma},F(S)\db{\sigma}} \circ Q(a^{F(S)}_{\Gamma,x:\sigma})(F(S)\db{\psi [\Gamma,x:\sigma]}) \\
		= \ & Q(a_{I}^{-1}) \circ \Omega_{F^{o}(I),X} \circ Q(a_{I} \times id_{F(X)}) \circ Q(a^{F(S)}_{\Gamma,x:\sigma})(F(S)\db{\psi [\Gamma,x:\sigma]}) \\
		= \ & Q(a_{I}^{-1}) \circ \Omega_{F^{o}(I),X} \circ Q(a_{I} \times id_{F(X)}) \circ Q(a^{F(S)}_{\Gamma,x:\sigma}) \circ Q(a_{S\db{\Gamma,x:\sigma}}^{-1})\circ F^{p}_{S \db{\Gamma,x:\sigma}}(S \db{\psi  [\Gamma,x:\sigma]}) \\
		= \ & Q(a_{I}^{-1}) \circ \Omega_{F^{o}(I),X} \circ Q(a_{I,X}^{-1}) \circ Q(F^{o}(a^{S}_{\Gamma,x:\sigma})) \circ F^{p}_{S \db{\Gamma,x:\sigma}}(S \db{\psi [\Gamma,x:\sigma]}) \\
		= \ & Q(a_{I}^{-1}) \circ \Omega_{F^{o}(I),X} \circ Q(a_{I,X}^{-1}) \circ F^{p}_{I \times X} \circ P(a^{S}_{\Gamma,x:\sigma})(S \db{\psi [\Gamma,x:\sigma]})  \\
		= \ & Q(a_{I}^{-1}) \circ F^{p}_{I} \circ \Omega_{I,X} \circ P(a^{S}_{\Gamma,x:\sigma})(S \db{\psi [\Gamma,x:\sigma]})   \\
		= \ & Q(a_{S \db{\Gamma}}^{-1}) \circ F^{p}_{S\db{\Gamma}}(S \db{\Omega_{x:\sigma}(\psi)[\Gamma]}. \qedhere
		\end{align*}
	\end{proof}
	
	The next result says we can identify $ T $-models $ S \in (\m{C},P) $ with morphisms $ \overline{S} \colon (\m{C}_{T},P_{T}) \to (\m{C},P) $ which we use to develop an ``algebraic'' view of entailment. 
	\begin{thrm}\label{thrm: Each T-model is up to iso a unique image of the generic T-model.}
		Let $ \m{L} $ be adequate and $ T $ an $ \m{L} $-theory. For each $ (\m{C},P) \in \Ob(\FA) $ and $ T $-model $ S $ in $ (\m{C},P) $ there is a morphism $ \overline{S} \colon (\m{C}_{T},P_{T}) \to (\m{C},P) $ in $ \FA $, unique up to a $ 2 $-isomorphism, such that $ \overline{S}(G) = S $, where $ G $ is the generic $ T $-model in $ (\m{C}_{T},P_{T}) $.
	\end{thrm}
	\begin{proof}
		We define $ \overline{S}^{o} $ on objects by $ \overline{S}^{o}(\Gamma) \coloneqq S \db{\Gamma} $, and on morphisms by 
		$$ \overline{S}^{o} ([M_{1}[\Gamma]/ \! \sim,\ldots, M_{n}[\Gamma] / \! \sim ]) \coloneqq \langle S \db{M_{1}[\Gamma]}, \ldots, S\db{M_{n}[\Gamma]} \rangle. $$
		Since $ S $ is a $ T $-model, $ \overline{S}^{o} $ is well defined on morphisms and it is straightforward to show that $ \overline{S}^{o} \colon \m{C}_{T} \to \m{C} $ is a product preserving functor such that $ \overline{S}(G) = S $, where $ G $ is the generic $ T $-algebra in $ \m{C}_{T} $. We extend $ \overline{S}^{o} $ by defining $ \overline{S}^{p}_{\Gamma}(\phi [\Gamma]/ \! \sim) \coloneqq S\db{\phi [\Gamma]} $, and $ \overline{S}^{p}_{\Gamma} $ is well-defined and monotone because $ S $ is a $ T $-model. By definition, $ \overline{S}(G) $ agrees with $ S $ on relation symbols. It is also straightforward to verify that $ \overline{S}^{p} \colon P_{T} \to P\circ \overline{S}^{o} $ is a natural transformation, and that for each context $ \Gamma $, $ \overline{S}^{P}_{\Gamma} $ is an $ \s{L} $-homomorphism.
		
		We now verify Condition 2. For $ m > 0 $, let $ \Gamma = x_{1}:\sigma_{1}, \ldots, x_{n}:\sigma_{n} $, $ \Gamma^{\prime} = x_{n+1}:\sigma_{n+1}, \ldots, x_{n+m}:\sigma_{n+m} $, and for each $ i \in \{1, \ldots, n+m \} $, let $ \pi_{i} \colon S \db{\Gamma,\Gamma^{\prime}} \to S \db{\sigma_{i}} $ be the $ i $th projection map. Let $ S \db{\Gamma} = I $, $ S \db{\Gamma^{\prime}} = I^{\prime} $ and for $ i \in \{ 1, \ldots, n+m \} $, let $ S \db{\sigma_{i}} = X_{i} $. Let $ a \colon S\db{\Gamma, \Gamma^{\prime}} \to S \db{\Gamma} \times S \db{\Gamma^{\prime}} $ be the change-of-product isomorphism. Then for each $ \Omega \in \s{L}_{q} $, if $ m > 0 $,
		\begin{align*}
		& \overline{S}^{p}_{\Gamma} \circ \Omega_{\Gamma,\Gamma^{\prime}}(\phi [\Gamma, \Gamma^{\prime}]/ \! \sim) \\
		= \ & \overline{S}^{p}_{\Gamma} ( \Omega_{ x_{n+1}:\sigma_{n+1}} \cdots \Omega_{x_{n+m}: \sigma_{n+m}}(\phi) [\Gamma]/ \! \sim) \\
		= \ & \Omega_{I,X_{n+1}} \circ \ldots \circ \Omega_{I \times X_{n+1} \ldots \times X_{n+m-1}, X_{n+m}}(S \db{\phi [\Gamma,\Gamma^{\prime}]} ) \\
		= \ & \Omega_{S \db{\Gamma}, S \db{\Gamma^{\prime}}} \circ a^{-1 *} ( S \db{\phi [\Gamma, \Gamma^{\prime}]}) \tag*{(by \ref{con: 5 of prop-cat})} \\
		= \ & \Omega_{\overline{S}^{o}(\Gamma), \overline{S}^{o}(\Gamma^{\prime})} \circ a^{-1 *} \circ \overline{S}^{p}_{\Gamma \times \Gamma^{\prime}}(\phi[\Gamma,\Gamma^{\prime}]/ \! \sim ).
		\end{align*}
		Otherwise, $ \Gamma^{\prime} = [ \, ] $ and
		\begin{align*}
		\overline{S}^{p}_{\Gamma} \circ \Omega_{\Gamma,[ \, ]}(\phi \, [\Gamma]/ \! \sim) & = \overline{S}_{\Gamma}^{p}(\phi \, [\Gamma]/ \! \sim) = S \db{\phi \, [\Gamma]} \\
		& = \Omega_{S \db{\Gamma},S\db{ \, }}\circ P(\pi_{1}^{S\db{\Gamma},S \db{ \, }}) ( S \db{\phi [\Gamma]}) \tag{by \ref{con: 5 of prop-cat} } \\
		& = \Omega_{\overline{S}^{o}\Gamma, \overline{S}^{o}[ \, ]} \circ P(a_{\overline{S},\Gamma,[ \, ]}^{-1}) \circ S^{p}_{\Gamma \times [ \, ]}(\phi \, [\Gamma]/ \! \sim).
		\end{align*}
		
		Let $ \Gamma = x_{1}:\sigma_{1}, \ldots, x_{n}:\sigma_{n} $, $ \Gamma^{\prime} = x_{n+1}:\sigma_{1}, \ldots, x_{2n}:\sigma_{n} $. If $ n > 0 $,
		\begin{align*}
		& \overline{S}^{p}_{\Gamma \times \Gamma^{\prime}}(Eq_{\Gamma}) \\
		= \ & \bigcon_{i=1}^{n} S\db{x_{i}=_{\sigma_{i}}x_{n+i}[\Gamma,\Gamma^{\prime}]} \\
		= \ & \bigcon_{i=1}^{n}P(\langle \pi_{i}^{S\db{\Gamma,\Gamma^{\prime}}},\pi_{i+n}^{S \db{\Gamma,\Gamma^{\prime}}}\rangle) Eq_{S\db{\sigma_{i}}} \\
		= \ & \bigcon_{i=1}^{n} P(\langle \pi_{i}^{S\db{\Gamma}}\pi_{1}^{S\db{\Gamma},S\db{\Gamma^{\prime}}}, \pi_{i}^{S\db{\Gamma^{\prime}}}\pi_{2}^{S\db{\Gamma},S\db{\Gamma^{\prime}}}\rangle \circ a )Eq_{S\db{\sigma_{i}}} \\
		= \ & P(a) \circ \bigcon_{i=1}^{n}  P(\langle \pi_{i}^{S\db{\Gamma}}\pi_{1}^{S\db{\Gamma},S\db{\Gamma^{\prime}}}, \pi_{i}^{S\db{\Gamma^{\prime}}}\pi_{2}^{S\db{\Gamma},S\db{\Gamma^{\prime}}}\rangle)Eq_{S\db{\sigma_{i}}}      \\
		= \ & P(a) \circ Eq_{\overline{S}^{o}(\Gamma)}. \tag*{(by \ref{con: 6 of prop-cat.})}
		\end{align*}
		And if $ n = 0 $,
		\begin{align*}
		\overline{S}^{p}_{[ \, ]}(Eq_{[ \, ]}) & = \overline{S}^{p}_{[ \, ]}(e[ \, ]/ \! \sim) =  S\db{e[ \, ]} = e_{1} \tag{ by \ref{con: 6 of prop-cat.} }  \\
		& = P(a_{\overline{S},[ \, ],[ \, ]})(e_{1 \times 1}) = P(a_{\overline{S},[ \, ],[ \, ]})(Eq_{\overline{S}^{o}[ \, ]}).
		\end{align*}
		
		Now all that remains to show is that if $ F \colon (\m{C}_T,P_{T}) \rightarrow (\m{C},P) $ is another morphism such that $ F(G) = \overline{S}(G) = S $ then they are 2-isomorphic.
		
		Let $ a_{F,\Gamma} \colon F^{o}(\Gamma) \to \overline{F(G)}^{o}(\Gamma) = \overline{S}^{o}(\Gamma) $ be the change-of-product isomorphism. Then $ a_{F} \colon F^{o} \Rightarrow \overline{S}^{o}$ is a natural isomorphism, and $ Pa_{F} \cdot \overline{S}^p = Pa_{F} \cdot \overline{F(G)}^{p} = F^{p} $, from Equation \ref{eqn: goal to prove morphisms preserve modesl}. Thus $ a_{F} $ is a $2$-isomorphism from $ F $ to $ \overline{S} $. 
	\end{proof}
	
	\begin{rem}
		From Theorem 5, morphisms between classifying prop-categories correspond to interpretations of one theory in another. Consider the following example: Let $ \m{L} $ be classical first-order logic, $ Sg $ the signature of rings with sort $ \sigma $, nullary operations $ 0,1 :\sigma $ and binary operations $ +,-,\cdot : \sigma, \sigma \to \sigma $. Let $ T $ be the $ Sg $-theory of fields with no square root of $ -1 $ and $ T^{\prime} $ the theory of fields with a square root of $ -1 $. Then one can define a $ T^{\prime} $-model $ S $ in $ (\m{C}_{T},P_{T}) $, such that $ S \db{\sigma} = [a:\sigma, b:\sigma] $, $$ S \db{0} = [ 0: \sigma \, [ \, ], 0 :\sigma \, [ \, ]], \quad S \db{1} = [1:\sigma \, [ \, ], 0 :\sigma \, [\,]],$$ $$ S \db{\pm} = [a \pm c : \sigma \, [\Gamma], b \pm d :\sigma \, [\Gamma]], \ \ \text{and} \ \ S \db{\cdot} = [ac-bd:\sigma \, [\Gamma], ad + bc :\sigma \, [\Gamma]], $$ where $ \Gamma = [a:\sigma,b:\sigma,c:\sigma,d:\sigma] $. This induces a morphism $ \overline{S} \colon (\m{C}_{T^{\prime}},P_{T^{\prime}}) \to (\m{C}_{T},P_{T}) $ which encodes the interpretation of $ T^{\prime} $ in $ T $ and gives, for each $ T $-model viewed as a morphism $ F \colon (\m{C}_{T},P_{T}) \to (\m{C},P) $ a $ T^{\prime} $-model $ F \circ \overline{S} \colon (\m{C}_{T^{\prime}},P_{T^{\prime}}) \to (\m{C}, P) $.
		
		There are more complicated notions of the interpretation of one theory in another which, given our current setup, cannot be captured by morphisms between classifying prop-categories. For a simple example, constructing the field of fractions out of an integral domain induces an interpretation of the theory of fields in the theory of domains. However, the ``set'' of field elements is a quotient of a subset of the pairs of domain elements. In order to capture this example as a prop-category morphism between the corresponding prop-categories one could extend the type theory to include subset and quotient types and require the corresponding additional structure on the prop-categories as done in \cite[p.\ 272]{Jacobs1999}.
	\end{rem}
	
	\begin{defn}\label{defn: Kernel of a morphism.}
		Let $ F \colon (\m{C},P) \to (\m{D},Q) $ be a morphism in $ \FA $. We define the \textbf{kernel of $ F $}, denoted $ \ker F $ to consist of the following data:
		\begin{enumerate}
			\item A relation on $ \Ob(\m{C}) $, such that $ c_{1} \sim c_{2} $ iff $ F^{o}c_{1} = F^{o}c_{2} $.
			\item A relation on $ \Mor(\m{C}) $ such that $ f_{1} \sim f_{2} $ iff $ F^{o}f_{1} = F^{o}f_{2} $.
			\item A relation on $ \bigsqcup_{c \in \Ob(\m{C})}P(c) $ such that $ r_{1} \prec r_{2} $ iff $ F^{p}_{c_{1}}(r_{1}) \leq F^{p}_{c_{2}}(r_{2}) $, where $ r_{i} \in P(c_{i}) $.
		\end{enumerate}
		If $ K \colon (\m{C},P) \to (\m{E},R) $ is another morphism, we say $ \ker K \leq \ker F $, to assert the relations of $ \ker K $ are contained in the corresponding relations of $ \ker F $.
	\end{defn}

	Let $ S $ be an $ Sg $-structure in $ (\m{C},P) \in \Ob(\FA) $ and let $ a \in \A_{Sg} $. If $ a = \phi_{1}, \ldots, \phi_{n} \vdash \phi \ [\Gamma] $, define
	$$
	S\db{a} \coloneqq (\bigotimes_{i=1}^{n}S\db{\phi_{i}[\Gamma]},S\db{\phi [\Gamma]}),
	$$
	and if $ a $ is an equation $ M_{1} = M_{2}: \tau \ [\Gamma] $, define
	$$
	S\db{a} \coloneqq (S \db{M_{1}: \tau \ [\Gamma]},S \db{M_{2}:\tau \ [\Gamma]}).
	$$
	For each $Sg$-theory $T$, we define $ S \db{T} \coloneqq \{ S\db{a}: a \in \A(T) \} $ and for $ F \colon (\m{C}_{T},P_{T}) \to (\m{C},P) $ where $ T $ is an $ \m{L}^{m} $-theory, let $ \Th(F) $ be the $ T $-theory such that
	$$
	\A(\Th(F)) \coloneqq \{ a \in \A_{Sg} : G\db{a} \in \ker F \}.
	$$
	It is straightforward to verify that $ \Th(F) = \Th(F(G)) $. Thus, for each structure $ S  \in (\m{C},P) $, $ \Th(S) = \Th(\overline{S}(G)) =  \Th(\overline{S}) $. Moreover, whenever $ F $ is $ 2 $-isomorphic to a parallel $ 2 $-cell $ K $, then $ \Th(F) = \Th(K) $.

	Let $ T $, and $ T^{\prime} $ be $ Sg $-theories and $ (\m{C},P) \in \Ob(\FA) $. From our prior observations the following are equivalent:
	\begin{enumerate}
		\item $ T \vDash_{(\m{C},P)} T^{\prime} $.
		\item $ \forall F \colon (\m{C}_{Sg},P_{Sg}) \to (\m{C},P) $ such that $ T \leq \Th(F) $, then $ T^{\prime} \leq \Th(F) $.
		\item $ \forall F \colon (\m{C}_{Sg},P_{Sg}) \to (\m{C},P) $ such that $ G\db{T} \subseteq \ker F $, then $ G \db{T^{\prime}} \subseteq \ker F $.
	\end{enumerate}
	Thus in the sequel, we will identify $ Sg $-structures in $ (\m{C},P) $ with morphisms $ F \colon (\m{C}_{Sg},P_{Sg}) \to (\m{C},P) $.
		
	\begin{rem}
	Prop-categorical semantics provide a natural notion of structural action on theories: Let $ H \colon (\m{C}_{Sg},P_{Sg}) \to (C_{Sg^{\prime}},P_{Sg^{\prime}}) $ be a morphism and let $ T_{1},T_{2} $ be $ Sg $-theories. We can define an action of $ H $ on $ Sg $ equations and sequents in context which can be extended to $ T_{1} $ by taking the union: that is $ H \cdot T_{1} $ is the $ Sg^{\prime} $ theory whose assertions are $ \bigcup_{a \in A(T_{1})} H \cdot a $. For equations-in-context $ M^{1} = M^{2} : \tau \ [\Gamma] $,
	$$
	H \cdot (M^{1} = M^{2}: \tau \ [\Gamma]) \coloneqq \{ M^{1}_{i} = M^{2}_{i}: \tau_{i} \ [\Gamma_{H}]  \}_{i=1}^{n},
	$$
	where $ H^{o}(M^{i}:\tau \ [\Gamma]) = [M^{i}_{1}:\tau_{1}\, [\Gamma_{H}],\ldots, M^{i}_{n}:\tau_{n} \, [\Gamma_{H}]] $ for $ i \in \{1,2 \} $. And for sequents-in-context $ \phi_{1}, \ldots, \phi_{n} \vdash \phi_{n+1} \ [\Gamma] $,
	$$
	H \cdot (\phi_{1}, \ldots, \phi_{n} \vdash \phi_{n+1} \ [\Gamma]) = \phi_{1}^{H}, \ldots, \phi_{n}^{H} \vdash \phi_{n+1}^{H} \ [\Gamma_{H}],
	$$
	where $ H^{p}_{\Gamma}(\phi_{i} \, [\Gamma]) = \phi_{i}^{H} \, [\Gamma_{H}] $, for all $ i \in \{ 1, \ldots, n+1 \} $. Since products in $ (\m{C}_{Sg},P_{Sg}) $ are unique up to permutation of the list of their variable sort pairs, if $ K \colon (\m{C}_{Sg^{\prime}},P_{Sg^{\prime}}) \to (\m{C}_{Sg^{\prime \prime}},P_{Sg^{\prime \prime}}) $, then $ (K \circ H) \cdot T_{1} = K \cdot (H \cdot T_{1}) $ and clearly, $ id_{(\m{C}_{Sg},P_{Sg})} \cdot T_{1} = T_{1} $. Moreover, suppose $ \m{V} \subseteq \Ob(\FA_{\s{L}}) $ and that $ T_{1} \vDash_{\m{V}} T_{2} $. Let $ (\m{C},P) \in \m{V} $ and $ F \colon (\m{C}_{Sg^{\prime}},P_{Sg^{\prime}}) \to (\m{C},P) $ such that $ G \db{ H \cdot T_{1} } \subseteq \ker F $. Then $ G \db{ T_{1} } \subseteq \ker F \circ H $ and since $ T_{1}\vDash_{\m{V}} T_{2} $, $ G \db{T_{2}} \subseteq \ker F \circ H $. It follows that $ G \db{H \cdot T_{2}} \subseteq \ker F $ and so $ H \cdot T_{1} \vDash_{\m{V}} H \cdot T_{2} $. Therefore, all logics defined semantically by subcollections of $ \Ob(\FA_{\s{L}}) $ are structural with respect to the actions of morphisms, a key ingredient to the main theory of AAL.
	\end{rem}

	\section{Basic Constructions in $ \FA $} \label{sec: Basic Constructions}
	
	Let $ (\m{C},P) \in \Ob(\FA) $ and consider the signature $ Sg $ whose sorts are $ \Ob(\m{C}) $, whose function symbols are $ f:c_{1}, \ldots,c_{n} \to c $ for each $ f: c_{1} \times \cdots \times c_{n} \to c \in \Mor(\m{C}) $, and whose relation symbols are $ R \subseteq c_{1}, \ldots,c_{n} $ for each $ R \in P(c_{1} \times \cdots \times c_{n}) $. Note that $ f : c_{1} \times \cdots \times c_{n} \to c $ is included both as an $ n $-ary operation symbol $ f : c_{1}, \ldots, c_{n} \to c $ and as a unary operation symbol $ f : c_{1} \times \cdots \times c_{n} \to c $ in $ Sg $ and similarly for relation symbols. There is a canonical $ Sg $-structure $ S $ in $ (\m{C},P) $ called the \textbf{internal structure of $ (\m{C},P) $}, where $ \db{c} = c $, $ \db{f} = f $ and $ \db{R} = R $ for each sort, function symbol and relation symbol respectively. We define $ \Th(\m{C},P) \coloneqq \Th(S) $.
	
	\begin{prop}\label{prop: 2-equivalence of a prop-category and its syntactic analogue}
		Let $ (\m{C},P) \in \Ob(\FA) $ and let $ S $ be the internal structure of $ (\m{C},P) $, and $ T $ the theory of $ S $. Then $ \overline{S} \colon (\m{C}_{T},P_{T}) \to (\m{C},P) $ determines a $ 2 $-equivalence.
	\end{prop}
	\begin{proof}
		Define $ \iota \colon (\m{C},P) \to (\m{C}_{T},P_{T}) $ by
		$$
		\iota^{o}(c) \coloneqq x:c, \ \ \iota^{o}(f) \coloneqq f(x):c_{2} \ [x:c_{1}]/ \! \sim, \ \text{ and } \ \iota^{p}(R) \coloneqq R(x) \ [x:c] / \! \sim,
		$$
		for $ c,c_{1},c_{2} \in \Ob(\m{C}) $, $ f \in \m{C}(c_{1},c_{2}) $ and $ R \in P(c) $. It is straightforward to verify that $ \iota^{o} $ is a finite product preserving functor. Let $ f \colon c_{1} \to c_{2} \in \Mor(\m{C}) $. Then for $ R \in P(c_{2}) $, 
		\begin{align*}
		\iota_{c_{1}}^{p} \circ P(f)(R) & = P(f)(R)(x)[x: c_{1}] / \! \sim \\
		& = R(f(x))[x:c_{1}]/ \! \sim \\
		& = P_{T}(f(x):c_{2}[x:c_{1}]/ \! \sim) (R(x)[x:c_{2}]/ \! \sim) \\
		& = P_{T}( \iota^{o}(f)) \circ \iota_{c_{2}}^{p}(R).
		\end{align*}
		Thus $ \iota^{p} $ is a natural transformation. It is also straightforward to verify for each $ c \in \Ob(\m{C}) $ that $ \iota^{p}_{c} $ is a monotone $ \s{L} $-algebra homomorphism.
		Let $ \Omega \in \s{L}_{q} $, $ c_{1},c_{2} \in \Ob(\m{C}) $ and $ R \in P(c_{1} \times c_{2}) $. Then
		\begin{align*}
		& \Omega_{\iota^{o}c_{1},\iota^{o}c_{2}} \circ P_{T}(a^{-1}) \circ \iota^{p}_{c_{1}\times c_{2}}(R) \\
		= \ & \Omega_{x_{2}:c_{2}}(R( \langle \pi_{1},\pi_{2} \rangle(x_{1},x_{2}) ))[x_{1}:c_{1}] / \! \sim \\
		= \ & \Omega_{c_{1},c_{2}}(R)(x)[x:c_{1}] / \! \sim \\
		= \ & \iota^{p}_{c_{1}} \circ \Omega_{c_{1},c_{2}}(R). 
		\end{align*}
		Let $ c \in \Ob(\m{C}) $. Then
		\begin{align*}
		& P_{T}(a^{-1}) \circ \iota^{p}_{c \times c}(Eq_{c}) \\
		= \ & Eq_{c}(\langle \pi_{1},\pi_{2} \rangle(x_{1},x_{2}))[x_{1}:c,x_{2}:c] / \! \sim \\
		= \ & Eq_{c}(x_{1},x_{2})[x_{1}:c,x_{2}:c] / \! \sim \\
		= \ & x_{1} =_{c} x_{2} \ [x_{1}:c,x_{2}:c] / \! \sim \\
		= \ & Eq_{x:c}.
		\end{align*}
		It follows that $ \iota $ is a morphism in $ \FA $.
		
		Now $ \overline{S}^{o} \circ \iota^{o} = id_{\m{C}} $, and for each $ c \in \Ob(\m{C}) $ and $ R \in P(c) $,
		$$
		(\overline{S}_{\iota^{o}}^{p} \cdot \iota^{p})_{c}(R) = \overline{S}_{\iota^{o}c}^{p} \circ \iota_{c}^{p}(R) = \overline{S}_{\iota^{o}c}^{p}(R(x)[x:c]/ \! \sim) = R.
		$$
		It follows that $ \overline{S} \circ \iota = \mathrm{id}_{(\m{C},P)} $.
		
		In the other direction, $$
		\iota^{o} \circ \overline{S}^{o}([x_{1}:c_{1}, \ldots, x_{n}:c_{n}]) = [x: c_{1} \times \cdots \times c_{n}],
		$$
		\begin{align*}
		& \iota^{o} \circ \overline{S}^{o}([M_{1}:c_{1}[\Gamma]/ \! \sim, \ldots, M_{m}:c_{m}[\Gamma] / \! \sim]) \\
		= \ & \langle S \db{M_{1}:c_{1}[\Gamma]}, \ldots S \db{M_{m}:c_{m}[\Gamma]}\rangle(x):c_{1} \times \cdots \times c_{m} \ [x:S \db{\Gamma}] / \! \sim
		\end{align*}
		and
		$$
		(\iota_{\overline{S}^{o}}^{p} \cdot \overline{S}^{p})_{\Gamma}(\phi [\Gamma]/ \! \sim) = S \db{\phi [\Gamma]}(x) [x: S\db{\Gamma}]/ \! \sim.
		$$
		For each context $ \Gamma = x_{1}:c_{1}, \ldots, x_{n}:c_{n} $, define $ \eta \colon \iota^{o} \cdot \overline{S}^{o} \Rightarrow \mathrm{id}_{C_{T}} $ by 
		$$ \eta_{\Gamma} = [\pi_{1}(x):c_{1}[x: S \db{\Gamma}]/ \! \sim, \ldots, \pi_{n}(x):c_{n}[x : S \db{\Gamma}]/ \! \sim] . $$
		One may verify that $ \eta $ is a natural isomorphism and  that $ P_{T} \eta \cdot \mathrm{id}^{p}_{(\m{C}_{T}, P_{T})} = \iota^{p}_{\overline{S}^{o}} \cdot \overline{S}^{p} $. It follows that $ \eta \colon \iota \circ \overline{S} \Rightarrow \mathrm{id}_{(\m{C}_{T},P_{T})} $ is a $ 2 $-isomorphism and so $ (\m{C}_{T},P_{T}) $ and $ (\m{C},P) $ are equivalent.
	\end{proof}
	
	Thus $ (\m{C}_{T},P_{T}) \equiv (\m{C},P) $ and so $ (\m{C}_{T},P_{T}) $ is a syntactic representation of $ (\m{C},P) $ similar to the representation of an algebra as the free algebra over its elements quotiented out by is equational theory. Moreover, it can be shown that if $ (\m{C},P),(\m{D},Q) \in \Ob(\FA) $ are equivalent in $ \FA $, then $ \vDash_{(\m{C},P)} \ = \ \vDash_{(\m{D},Q)} $ and so $ \vDash_{(\m{C},P)} \ = \ \vDash_{(\m{C}_{T},P_{T})} $.
	
	Let $ (\m{C},P),(\m{D},Q) \in \Ob(\FA) $. We say \textbf{$ (\m{D},Q) $ is a sub-prop-category of $ (\m{C},P) $}, if there exists a morphism $ \iota \colon (\m{D},Q) \to (\m{C},P) $, such that $ \iota^{o} $ is faithful and for each $ c \in \Ob(\m{C}) $, $ \iota^{p}_{c} $ is an order embedding. We call the morphism $ \iota $ a \textbf{sub-prop-morphism}.

	\begin{prop}\label{prop: FAL is closed under subprop-categories}
		If $ (\m{D},Q) $ is a sub-prop-category of $ (\m{C},P) $, then $ \vDash_{(\m{C},P)} \, \subseteq \, \vDash_{(\m{D},Q)} $.
	\end{prop}
	\begin{proof}
		Suppose $ T \vDash_{(\m{C},P)} T^{\prime} $ and consider $ F \colon (\m{C}_{Sg},P_{Sg}) \to (\m{D},Q) $, such that $ G\db{T} \subseteq \ker F $. Then $ G \db{T} \subseteq \ker(\iota \circ F) $, and since $ T \vDash_{(\m{C},P)} T^{\prime} $, $ G\db{T^{\prime}} \subseteq \ker(\iota \circ F) $. Since $ \iota $ is a sub-prop-morphism,  $ G\db{T^{\prime}} \subseteq \ker(F) $ and so $ T \vDash_{(D,Q)} T^{\prime} $.
	\end{proof}

	Let $ F \colon (\m{C},P) \to (\m{D},Q) $ be a morphism in $ \FA $. In general, $ (F(\m{C}),F(P)) $, is not a sub-prop-category of $ (\m{D},Q) $ where $ F(\m{C}) $ is the image of $ \m{C} $ under $ F^{o} $, and for all $ g \colon d_{1} \to d_{2} \in \Mor(F(\m{C})) $, $ F(P)(g) \coloneqq Q(g)|_{F^{p}_{c_{2}}P(c_{2})} \colon F^{p}_{c_{2}}P(c_{2}) \to F^{p}_{c_{1}}P(c_{1}) $, where $ F^{o}c_{i} = d_{i} $ for $ i \in \{1, 2 \} $. When $ (F(\m{C}),F(P)) $ does define a sub-prop-category of $ (\m{D},Q) $ we call it the \textbf{image of $ F $}.
	
	\begin{lem}\label{prop: The image of an injective on objects morphism determines a sub-prop-category.}
		Let $ F \colon (\m{C},P) \to (\m{D},Q) $ be a morphism in $ \FA $ such that $ F^{o} $ is injective on objects. Then $ (F(\m{C}), F(P) ) $ is a sub-prop-category of $ (\m{D},Q) $ and there exists a unique morphism $ H \colon (\m{C},P) \to (F(\m{C}),F(P)) $ such that $ \iota \circ H = F $, where $ \iota \colon (F(\m{C}),F(P)) \hookrightarrow (\m{D},Q) $ is the inclusion morphism. We call $ H $ the \textbf{corestriction of $ F $ to $ (F(\m{C}),F(P)) $}. Furthermore, $ H $ is strictly finite product preserving, full, surjective on objects, and for each $ c \in \Ob(\m{C}) $, $ H^{p}_{c} $ is surjective. 
	\end{lem}
	\begin{proof}
		Since $ F^{o} $ is injective on objects, $ F(\m{C}) $ is a sub-category of $ \m{D} $. Let $ d_{1},d_{2} \in \Ob(F(\m{C})) $ and $ F^{o}(c_{i}) = d_{i} $. Then the designated product diagram of $ d_{1} $ and $ d_{2} $ in $ F(\m{C}) $, which we denote $ (d_{1} \hat{\times} d_{2}, \hat{\pi}_{1},\hat{\pi}_{2}) $ (to distinguish it from the designated diagram of $ d_{1} \times d_{2} $ in $ \m{D} $) is $ (F^{o}(c_{1} \times c_{2}), F^{o}(\pi_{1}),F^{o}(\pi_{2})) $ where $ (c_{1} \times c_{2}, \pi_{1}, \pi_{2}) $ is the designated product diagram of $ c_{1} \times c_{2} $ in $ \m{C} $. Given $ g_{1} \colon d \to d_{1} $ and $ g_{2} \colon d \to d_{2} $ we let $ < \! g_{1},g_{2} \! > \, \coloneqq  F^{o}(\langle f_{1},f_{2} \rangle) $, which is the unique morphism such that $ \hat{\pi}_{i} \circ < \! g_{1},g_{2} \! > \, = g_{i} $ for $ i \in \{ 1,2 \}. $ We also define the designated terminal object in $ F(\m{C}) $ to be $ 1_{F(\m{C})} \coloneqq F^{o}(1_{\m{C}}) $. Define $ F(P)(d_{1}) \coloneqq F^{p}(P(c_{1})) \leq Q(d_{1}) $ and for $ g \in F(\m{C})(d_{1},d_{2}) $, if $ F^{o}(f) = g $, define $ F(P)(g) \coloneqq Q(g)|_{F(P)(d_{2})} = Q(F^{o}f)|_{F^{p}_{c_{2}}P(c_{2})} $. Since $ F^{p} \colon P \Rightarrow Q\circ F^{o} $, we have $ F(P)(g)(F(P)(d_{2})) \subseteq F(P)(d_{1}) $. For $ \Omega \in \s{L}_{q} $ we denote the interpretation of $ \Omega $ in $ (F(\m{C}),F(P)) $ by $ \Omega^{F(P)} $, which we define for $ r \in F(P)(d_{1} \hat{\times} d_{2}) $ as 
		$$ \Omega^{F(P)}_{F^{o}c_{1},F^{o}c_{2}}(r) \coloneqq \Omega_{F^{o}c_{1},F^{o}c_{2}} \circ Q(a_{F,c_{1},c_{2}}^{-1})(r). $$
		Let $ r^{\prime} \in P(c_{1} \times c_{2}) $ such that $ F^{p}_{c_{1} \times c_{2}}(r^{\prime}) = r $. Then,
		$$
		\Omega^{F(P)}_{d_{1},d_{2}}(r) = \Omega_{d_{1},d_{2}}\circ Q(a^{-1})\circ F^{P}_{c_{1}\times c_{2}}(r) = F^{p}_{c_{1}}(\Omega_{c_{1},c_{2}}(r^{\prime})) \in F(P)(d_{1}).
		$$
		Equality $ Eq^{F(P)} $ is defined by 
		$$ Eq_{F^{o}c}^{F(P)} \coloneqq Q(a_{F,c,c}) (Eq_{F^{o}c}) =  F_{c \times c}^{p}(Eq_{c}) \in F(P)(F^{o}c \, \hat{\times} \, F^{o}c). $$
		
		First we show that $ (F(\m{C}),F(P)) \in \FA $. Conditions 1 and 2 are immediate. Condition 3 follows from the fact that for each $ F^{o}c \in \Ob(F(\m{C})) $, 
		$$ \Omega^{F(P)}_{(\cdot), d} = \Omega_{(\cdot),F^{o}c} \cdot Q(a_{F,(\cdot),c}^{-1})|_{F(P)((\cdot) \hat{\times} F^{o}c)}. $$ Condition 4 is satisfied as we take $ \otimes^{F(P)(d_{1})} \coloneqq \otimes^{Q(d_{1})}|_{F(P)(d_{1})} $ and $ e_{d_{1}}^{F(P)} \coloneqq F^{p}_{c_{1}}(e_{c_{1}}) = e_{d_{1}}  $. Next we consider Condition 5. Let $ d = F^{o}c \in F(\m{C}) $ and $ \Omega \in \s{L}_{q} $. Then 
		\begin{align*}
		\Omega_{d,F^{o}(1)}^{F(P)} \circ F(P)(\hat{\pi}_{1}^{d,F^{o}(1)}) & = \Omega_{d,F^{o}(1)} \circ Q(a_{F,c,1}^{-1}) \circ Q(\hat{\pi}_{1}^{d,F^{o}(1)})|_{F(P)(d)} \\
		\Omega_{d,F^{o}(1)} \circ Q(\pi_{1}^{d,F^{o}(1)})|_{F(P)(d)} & = id_{F(P)(d)}.
		\end{align*}
		Let $ d_{1},d_{2},d_{3} \in F(\m{C}) $, and $ F^{o}(c_{i}) = d_{i} $. Let $ r \in F(P)((d_{1} \hat{\times} d_{2})\hat{\times}d_{3}) $ and $ r^{\prime} \in P((c_{1} \times c_{2}) \times c_{3}) $ such that $ F^{p}_{(c_{1} \times c_{2}) \times c_{3}}(r^{\prime}) = r $. Then 
		\begin{align*}
		& = \Omega^{F(P)}_{d_{1}, d_{2}\hat{\times}d_{3}} \circ Q(\hat{a}_{d_{1}, d_{2},d_{3}})(r) \\
		& = \Omega_{d_{1}, d_{2}\hat{\times}d_{3}} \circ Q(a_{F,c_{1},c_{2}\times c_{3}}^{-1}) \circ Q(F^{o}(a_{c_{1},c_{2},c_{3}})) \circ F^{p}_{(c_{1} \times c_{2})\times c_{3}}(r^{\prime})  \\
		& = \Omega_{F^{o}c_{1},F^{o}(c_{2} \times c_{3})} \circ Q(a_{F,c_{1},c_{2} \times c_{3}}^{-1}) \circ F^{p}_{c_{1} \times (c_{2} \times c_{3})} \circ P(a_{c_{1},c_{2},c_{3}})(r^{\prime})\\
		& = F^{p}\circ \Omega_{c_{1},c_{2} \times c_{3}} \circ P(a_{c_{1},c_{2},c_{3}})(r^{\prime}) \\
		& = F^{p}_{c_{1}} \circ \Omega_{c_{1},c_{2}} \circ \Omega_{c_{1} \times c_{2},c_{3}}(r^{\prime}) \\
		& = \Omega_{F^{o}c_{1},F^{o}c_{2}}\circ Q(a^{-1}_{F,c_{1},c_{2}}) \circ F^{p}_{c_{1} \times c_{2}} \circ \Omega_{c_{1} \times c_{2},c_{3}}(r^{\prime}) \\
		& = \Omega_{F^{o}c_{1},F^{o}c_{2}} \circ Q(a_{F,c_{1},c_{2}}^{-1}) \circ \Omega_{F^{o}(c_{1} \times c_{2}),F^{o}c_{3}} \circ Q(a^{-1}_{F,c_{1} \times c_{2}, c_{3}}) \circ F^{p}_{(c_{1} \times c_{2})\times c_{3}}(r^{\prime}) \\
		& = \Omega_{d_{1},d_{2}}^{F(P)} \circ \Omega_{(d_{1} \hat{\times}d_{2}),d_{3}}^{F(P)}(r).
		\end{align*}
		Thus Condition 5 is satisfied. We now consider Condition 6:
		\begin{align*}
		Eq^{F(p)}_{F^{o}(1)} = F^{p}_{1 \times 1}(Eq_{1}) = F_{1 \times 1}^{p}(e_{1 \times 1}) = e_{F^{o}1 \hat{\times} F^{o}1}.
		\end{align*}
		Let $ d_{1},d_{2} \in \Ob(F(\m{C})) $ and $ F^{o}c_{i} = d_{i} $. Let $ c = c_{1} \times c_{2} $ and $ d = d_{1} \hat{\times} d_{2} $. Then,
		\begin{align*}
		Eq_{d}^{F(P)} & = F^{p}_{c \times c}(Eq_{c}) \\
		& =  F^{p}_{c \times c}(\bigcon_{i=1}^{2} P(\langle \pi_{i}^{c}\pi_{1}^{c,c}, \pi_{i}^{c}\pi_{2}^{c,c} \rangle)Eq_{c_{i}}) \\
		& = \bigcon_{i=1}^{2}  F^{p}_{c \times c}\circ P(\langle \pi_{i}^{c}\pi_{1}^{c,c}, \pi_{i}^{c}\pi_{2}^{c,c} \rangle)Eq_{c_{i}} \\
		& = \bigcon_{i=1}^{2}  Q( F^{o} \langle \pi_{i}^{c}\pi_{1}^{c,c}, \pi_{i}^{c}\pi_{2}^{c,c} \rangle )\circ F^{p}_{c_{i} \times c_{i}}Eq_{c_{i}} \\
		& = \bigcon_{i = 1}^{2} F(P)(< \! \hat{\pi}_{i}^{d} \hat{\pi}_{1}^{d,d},\hat{\pi}_{i}^{d} \hat{\pi}_{2}^{d,d} \! >) Eq_{d_{i}}^{F(P)}.
		\end{align*}
		It follows that $ (F(\m{C}),F(P)) \in \FA $.
		
		Let $ \iota^{o} \colon F(\m{C}) \hookrightarrow \m{D} $ be the inclusion functor and for each $ d \in \Ob(F(\m{C})) $, let $ \iota_{d}^{p} \colon F(P)(d) \hookrightarrow Q(d) $ be the inclusion homomorphism. It is immediate that $ \iota^{p} \colon F(P) \Rightarrow Q \circ \iota^{o} $ is natural. Let $ \Omega \in \s{L}_{q} $, $ d_{1},d_{2} \in \Ob(F(\m{C})) $ and $ F^{o}(c_{i}) = d_{i} $. For $ r \in F(P)(d_{1} \hat{\times}d_{2}) $,
		\begin{align*}
		\iota^{p}_{d_{1}} \circ \Omega_{d_{1},d_{2}}^{F(P)}(r) & = \Omega_{d_{1},d_{2}} \circ Q(a^{-1}_{F,c_{1},c_{2}})(r) \\
		& = \Omega_{\iota^{o}_{d_{1}},\iota^{o}_{d_{2}}} \circ Q(a^{-1}_{F,c_{1},c_{2}}) \circ \iota_{d_{1} \hat{\times}d_{2}}^{p}(r).
		\end{align*}
		And, 
		$$
		\iota_{d \hat{\times} d}^{p}(Eq^{F(P)}) = Q(a_{F,c,c}) \circ Eq_{\iota^{o}d}.
		$$
		It follows that $ \iota = (\iota^{o},\iota^{p}) \in \FA( (F(\m{C}),F(P)), (\m{D},Q)) $.
		
		Define $ H = (H^{o},H^{p}) $ by $ H^{o} \colon \m{C} \to F(\m{C}) $, where $ H^{o}(c) \coloneqq F^{o}(c) $, $ H^{o}(f) \coloneqq F^{o}(f) $ and $ H^{p}_{c} $ is the corestriction of $ F^{p}_{c} $ to $ F(P)(F^{o}c) $. The naturality of $ H^{p} \colon P \Rightarrow F(P)\circ H^{o} $ follows from the naturality of $ F^{p} \colon P \Rightarrow Q \circ F^{o} $. By definition $ H^{o} $ strictly preserves finite products and $ \iota \circ H = F $. Since $ \iota^{o} $ is injective on objects and faithful, and $ \iota^{p}_{d} $ is injective, for each $ d \in \Ob(F(\m{C})) $, $ H $ is uniquely defined so that $ H \circ \iota = F $. Let $ c_{1}, c_{2} \in \Ob(\m{C}) $ and $ r \in P(c_{1} \times c_{2}) $. That $ H $ commutes with the quantifiers and equality follows from the fact that $ F $ does. Thus, $ H \in \FA( (\m{C},P), (F(\m{C}),F(P))  ) $.
	\end{proof}

	\begin{prop}
		$ \FA $ has all products.
	\end{prop}
	
	\begin{proof}
		The terminal object $ 1_{\FA}$ is $ (*,1_{*})$ where $ * $ is the terminal category and $ 1_{*} $ is its unique endofunctor. It is straightforward to verify that $ (*,1_{*}) \in \FA $ and that for each prop-category $ (\m{C},P) \in \FA $, there is a unique morphism $ !_{(\m{C},P)} \colon (\m{C},P) \to (*,1_{*}) $, where $ !_{(\m{C},P)} = (!_{\m{C}},!_{P}) $.
		
		Let $ \{ (\m{C}_{i},P_{i}) \}_{i \in I} \subseteq \Ob(\FA) $ be non-empty. Let $ \prod_{I}\m{C}_{i} $ denote the product of the categories $ \{ \m{C}_{i} \}_{i \in I} $, that is, the category whose objects $ a $ specify for each $ i \in I $, an object $ a_{i} \in \Ob(\m{C}_{i}) $, and whose morphisms $ f \colon a \to b $ specify for each $ i \in I $ a morphism $ f_{i} \in \m{C}_{i}(a_{i},b_{i}) $. Composition $ g \circ f $ of morphisms $ f \colon a \to b $ and $ g \colon b \to c $ is defined for each $ i \in I $ as $ (g\circ f)_{i} \coloneqq g_{i} \circ f_{i} \colon a_{i} \to c_{i} $.  Since each $ \m{C}_{i} $ has finite products, we define $ 1 $ so that $ 1_{i} = 1_{\m{C}_{i}} $ and for $ a,b \in \Ob(\prod \m{C}_{i}) $, we define $ a \times b $ so that $ (a \times b)_{i} \coloneqq a_{i} \times b_{i} $. For $ k \in \{ 1,2 \} $, we define the projection maps $ \pi_{k}^{a,b} $ by $ \pi_{k,i}^{a,b} = \pi_{k}^{a_{i},b_{i}} $. One verifies that this determines a product diagram for $ a \times b $, where for every pair of maps $ f \colon c \to a $ and $ g \colon c \to b $, $ \langle f,g \rangle_{i} \coloneqq \langle f_{i},g_{i} \rangle $.
		
		We define $ \prod_{I}P_{i} \colon \prod_{I}\m{C}_{i} \to \Pos $, on objects $ a \in \Ob(\prod_{I}\m{C}_{i}) $ by $ (\prod_{I}P_{i})(a) \coloneqq \prod_{I}P_{i}(a_{i}) $ and on morphisms $ f \in \prod_{I}\m{C}_{i}(a,b) $ by $ (\prod_{I}P_{i})(f) \coloneqq \prod_{I}P_{i}(f_{i}) $. Then for $ f \colon a \to b$ and $ g \colon b \to c $, 
		$$ (\prod_{I}P_{i})(g\circ f) \coloneqq \prod_{I}P_{i}(f_{i})P_{i}(g_{i}) \coloneqq \prod_{I}P_{i}(f_{i}) \prod_{I}P_{i}(g_{i}) = (\prod_{I}P_{i})(f) \circ (\prod_{I}P_{i})(g).$$
		One also verifies that $ \prod_{I}P_{i}(id_{a}) = id_{\prod_{I}P_{i}(a)} $, so that $ \prod_{I}P_{i} $ is a contravariant functor and thus $ (\prod_{I}\m{C}_{i},\prod_{I}P_{i}) $ is a prop-category.
		
		For each $ a \in \Ob(\prod_{I}\m{C}_{i}) $, we define $ Eq_{a} $ and $ e_{a} $ so that $ Eq_{a,i} \coloneqq Eq_{a_{i}} $ and $ e_{a,i} \coloneqq e_{a_{i}} $. Let $ b \in \Ob(\prod_{I}\m{C}_{i}) $ and $ c = a \times b $. Then $ Eq_{1} = e_{1 \times 1} $ since for each $ i \in I $, 
		$$ Eq_{1,i} = Eq_{1_{\m{C}_{i}}} = e_{1_{\m{C}_{i}} \times 1_{\m{C}_{i}}} = e_{1 \times 1,i}. $$
		Also,
		\begin{equation*}
		Eq_{c} = (\prod_{I}P_{i})\langle \pi_{1}^{a,b}\pi_{1}^{c,c}, \pi_{1}^{a,b}\pi_{1}^{c,c} \rangle Eq_{a} \con (\prod_{I}P_{i})\langle \pi_{2}^{a,b}\pi_{1}^{c,c}, \pi_{2}^{a,b}\pi_{1}^{c,c} \rangle Eq_{b},
		\end{equation*}
		since for each $ j \in I $, 
		\begin{align*}
		& ( (\prod_{I}P_{i})\langle \pi_{1}^{a,b}\pi_{1}^{c,c}, \pi_{1}^{a,b}\pi_{1}^{c,c} \rangle Eq_{a} \con (\prod_{I}P_{i})\langle \pi_{2}^{a,b}\pi_{1}^{c,c}, \pi_{2}^{a,b}\pi_{1}^{c,c} \rangle Eq_{b} )_{j} \\
		= \ & P_{j}\langle \pi_{1}^{a_{j},b_{j}}\pi_{1}^{c_{j},c_{j}}, \pi_{1}^{a_{j},b_{j}}\pi_{1}^{c_{j},c_{j}} \rangle Eq_{a_{j}} \con P_{j} \langle \pi_{2}^{a_{j},b_{j}}\pi_{1}^{c_{j},c_{j}}, \pi_{2}^{a_{j},b_{j}}\pi_{1}^{c_{j},c_{j}} \rangle Eq_{b_{j}} \\
		= \ & Eq_{c,j}
		\end{align*}
		
		For $ \Omega \in \s{L}_{q} $, $ \Omega_{a,b} \colon \prod_{I}P_{i}(a \times b) \to \prod_{I}P_{i}(a) $ is defined so that $ \Omega_{a,b,i} \coloneqq \Omega_{a_{i},b_{i}} $. As with the conditions for equality, the remaining conditions for the quantifiers follow from the fact that they hold coordinate-wise. Therefore, $ (\prod_{I}\m{C}_{i},\prod_{I}P_{i}) \in \Ob(\FA) $.
		
		For $ j \in I $, the projection morphism $ \pi_{j} \colon (\prod_{I}\m{C}_{i}, \prod_{I}P_{i} ) \to (\m{C}_{j},P_{j}) $ consists of the projection functor $ \pi_{j}^{o} \colon \prod_{I}\m{C}_{i} \to \m{C}_{j} $ and the natural transformation $ \pi_{j}^{p} \colon \prod_{I}P_{i} \Rightarrow P_{j} \circ \pi_{j}^{o} $ defined for each object $ a \in \Ob(\prod_{I}\m{C}_{i}) $, by $ \pi_{j,a}^{p} \colon \prod_{I}P_{i}(a) \to P_{j}(a_{j}) $, i.e.\ $ \pi_{j,a}^{p} $ is the (monotone) projection homomorphism. For all $ f \colon a \to b \in \Mor(\prod_{I}\m{C}_{i}) $,
		\begin{equation*}
		\begin{tikzcd}[column sep = 40, row sep = 40]
		b \ar{d}[swap]{ \textstyle f^{op}} & \prod_{I}P_{i}(b_{i}) \ar{d}[swap]{\textstyle \prod_{I}P_{i}(f_{i}) } \ar{r}{ \textstyle \pi^{p}_{j,b}} & P_{j}(b_{j}) \ar{d}[swap]{\textstyle P_{j}(f_{j})} \\
		a & \prod_{I}P_{i}(a_{i}) \ar{r}{ \textstyle \pi^{p}_{j,a}} & P_{j}(a_{j})
		\end{tikzcd}
		\end{equation*}
		commutes and so $ \pi^{p}_{j} \colon \prod_{i}P_{i} \to P_{j} \circ \pi^{o}_{j} $ is a natural transformation. For all $ a \in \Ob(\prod_{I}\m{C}_{i}) $, 
		$$
		\pi^{p}_{j,a \times a}(Eq_{a}) = Eq_{a,j} = Eq_{a_{j}} = Eq_{\pi^{o}_{j}(a)}.
		$$
		For $ \Omega \in \s{L}_{q} $, $ a,b \in \Ob(\prod_{I}\m{C}_{i}) $ and $ A \in (\prod_{I}P_{i})(a \times b) $,
		$$
		\pi^{p}_{j,a}\circ \Omega_{a,b}(A) = \Omega_{a_{j},b_{j}}(A_{j}) = \Omega_{\pi_{j}(a),\pi_{j}(b)} \circ \pi^{p}_{j,a \times b}(A).
		$$
		It follows that $ \pi_{j} \colon (\prod_{I}\m{C}_{i}, \prod_{I}P_{\m{C}_{i}}) \to (\m{C}_{j},P_{\m{C}_{j}}) $ is a prop-category morphism.
		
		Let $ (\m{D},Q) \in \FA $ and for each $ j \in I $, let $ F_{j} \colon (\m{D}, Q) \to (\m{C}_{j},P_{j}) $ be a morphism in $ \FA $. Define $ \langle F_{i} \rangle_{I} \colon (\m{D},Q) \to (\prod_{I}\m{C}_{i},\prod_{I}P_{i}) $, so that $ \langle F_{i} \rangle_{I}^{o} \coloneqq \langle F^{o}_{i} \rangle_{I} $ and define $ \langle F_{i} \rangle_{I}^{p} \colon Q \Rightarrow \prod_{I}P_{i} \circ \langle F_{i} \rangle^{o}_{I} $ for $ d \in \Ob(\m{D}) $, by $ \langle F_{i} \rangle^{p}_{I,d} \coloneqq \langle F_{i,d}^{p} \rangle_{I} $. Then for all $ g \colon d_{1} \to d_{2} \in \Mor(\m{D}) $,
		\begin{equation*}
		\begin{tikzcd}[column sep = 40, row sep = 40]
		d_{2} \ar{d}[swap]{ \textstyle g^{op}} & Q(d_{2}) \ar{d}[swap]{\textstyle Q(g) } \ar{r}{ \textstyle \langle F^{p}_{i,d_{2}} \rangle_{I}} & \prod_{I}(P_{i}\circ F_{i}^{o}(d_{2})) \ar{d}[swap]{\textstyle \prod_{I}(P_{i} \circ F_{i}^{o}(g)) } \\
		d_{1} & P_{\m{D}}(d_{1}) \ar{r}{ \textstyle \langle F^{p}_{i,d_{1}} \rangle_{I} } & \prod_{I}(P_{i}\circ F^{o}_{i}(d_{1}))
		\end{tikzcd}
		\end{equation*}
		commutes since
		\begin{align*}
		\prod_{I}(P_{i} \circ F_{i}^{o}(g)) \circ \langle F_{i,d_{2}}^{p} \rangle_{I} & = \langle P_{i} \circ F^{o}_{i}(g) \circ F^{p}_{i,d_{2}} \rangle_{I} \\
		& = \langle F_{i,d_{1}}^{p} \circ Q(g) \rangle_{I} \\
		& = \langle F_{i,d_{1}}^{p} \rangle_{I} \circ Q(g). 
		\end{align*}
		It follows that  $ \langle F_{i} \rangle_{I}^{p} \colon Q \Rightarrow \prod_{I}P_{i} \circ \langle F_{i} \rangle^{o}_{I} $ is a natural transformation.
		
		Let $ d \in \Ob(\m{D}) $ and $ a \in \Mor(\prod_{I}C_{i}) $ such that for each $ j \in I $, $ a_{j} \colon F^{o}_{j}(d \times d) \to F_{j}^{o}d \times F_{j}^{o}d $ is the change-of-product isomorphism. Then for each $ j  \in I$,
		$$
		\langle F_{i,d \times d}^{p} \rangle_{I}(Eq_{d})_{j} = F_{j,d \times d}^{p} Eq_{d}  = P_{j}(a_{j})Eq_{F^{o}_{j}d},
		$$
		and so 
		$$
		\langle F_{i} \rangle^{p}_{I,d \times d}Eq_{d} = \prod_{I}P_{i}(a) Eq_{\langle F_{i} \rangle^{o}_{I}d}
		$$
		
		Let $ d_{1},d_{2} \in \Ob(\m{D}) $ and $ a \in \Mor(\prod_{I}\m{C}_{i}) $ such that for each $ j \in I $, $ a_{j} \colon F^{o}_{j}(d_{1} \times d_{2}) \to F^{o}_{j}d_{1} \times F^{o}_{j}d_{2} $ is the change-of-product isomorphism. Then, for $ \Omega \in \s{L}_{q} $,
		\begin{align*}
		\langle F_{i} \rangle_{I,d_{1}}^{p} \circ \Omega_{d_{1},d_{2}} & = \langle F_{i,d_{1}}^{p} \circ \Omega_{d_{1},d_{2}} \rangle_{I} \\
		& = \langle \Omega_{F^{o}_{i}d_{1},F^{o}_{i}d_{2}} \circ P(a_{i}^{-1}) \circ F_{i,d_{1} \times d_{2}}^{p}\rangle_{I} \\
		& = \Omega_{\langle F_{i}\rangle^{o}_{I}d_{1},\langle F_{i}\rangle^{o}_{I}d_{2}} \circ \prod_{I}P_{i}(a^{-1}) \circ \langle F_{i} \rangle_{I,d_{1} \times d_{2}}^{p}.
		\end{align*}
		Therefore, $ \langle F_{i} \rangle_{I} $ is a morphism.  Also, for each $ j \in I $ and $ d \in \Ob(\m{D}) $,
		$$
		(\pi^{p}_{j, \langle F_{i} \rangle_{I}^{o}} \cdot \langle F_{i} \rangle_{I}^{p} )_{d} = \pi^{p}_{j,\langle F_{i} \rangle_{I}^{o}d} \circ \langle F_{i,d}^{p} \rangle_{I} = F_{j,d}^{p},
		$$
		and 
		$$
		\pi_{j}^{o} \circ \langle F_{i} \rangle_{I}^{o} = F_{j}^{o}. 
		$$
		Thus for each $ j \in I $, $ \pi_{j} \circ \langle F_{i} \rangle_{I} = F_{j} $. If $ \phi \colon (\mathcal{D},Q) \to  (\prod_{I}\m{C}_{I},\prod_{I}P_{i}) $ is another morphism such that for each $ j \in I $, $ \pi_{j} \circ \phi = F_{j} $, then $ \phi^{o} = F_{i}^{o} $ and for each $ d \in \Ob(\m{D}) $, 
		$$
		(\pi^{p}_{j,\phi^{o}} \cdot \phi^{p})_{d} = \pi^{p}_{j,\langle F^{o}_{i}\rangle_{I}d} \circ \phi_{d}^{p} = F_{j,d}^{p}.
		$$
		It follows that $ \phi^{p}_{d} = \langle F_{i} \rangle^{p}_{I,d} $, and so $ \phi = \langle F_{i} \rangle _{I} $.
	\end{proof}

	\begin{prop}\label{prop: Logic of the product is stronger than the logic of the set.}
		Let $ \s{L} $ be a first-order language and $ \{ (\m{C}_{i},P_{i}) \}_{i \in I} \subseteq \Ob(\FA) $. Then
		$$ \vDash_{(\prod_{I}\m{C}_{i},\prod_{I}P_{i})} \, \supseteq \, \vDash_{\{ (\m{C}_{i},P_{i}) \}_{i \in I}}. $$
	\end{prop}
	\begin{proof}
		Suppose $ T \vDash_{\{ (\m{C}_{i},P_{i}) \}_{i \in I} } T^{\prime} $ and let $ S \in (\prod_{I}\m{C}_{i},\prod_{I}P_{i}) $ be a $ T $-model. Then for each $ i \in I $, $ \pi_{i}(S) $ is a $ T $-model in $ (\m{C}_{i},P_{i}) $. By assumption, $ \pi_{i}(S) $ satisfies $ T^{\prime} $  for all $ i \in I $ and it follows that $ S $ satisfies $ T^{\prime} $.
	\end{proof}
	As a consequence, for each logic $ \m{L} $, $ \FA_{\m{L}} $ has all products.

	\section{Fibered Homomorphism Theorem}\label{sec: Fibered Homomorphism Theorem}
	In this section we present a fibered homomorphism theorem which closely mirrors the corresponding result from universal algebra. These results give an orthogonal factorization system for $ \FA_{\m{L}} $. For a general reference on factorization systems see \cite{Riehl2008}.
	
	Fix a logic $ \m{L} $ and let $ \s{E} $ be the class of morphisms $ \epsilon $ in $ \FA_{\m{L}} $ such that $ \epsilon^{o} $ is bijective on objects, full and for each $ c \in \Ob(\m{C}) $, $ \epsilon^{p}_{c} $ is surjective. Let $ \s{M} $ be the collection of sub-prop-morphisms in $ \FA_{\m{L}} $.
	
	\begin{thrm}[Fibered Homomorphism Theorem]\label{thm: Fibered homomorphism theorem}
		Let $ F \colon (\m{C},P) \to (\m{D},Q) $ be a morphism in $ \FA_{\m{L}} $.
		\begin{enumerate}
			\item $ F $ factors as $ \psi \circ \epsilon $, for some $ \epsilon \in \s{E} $ and $ \psi \in \s{M} $.
			\item If $ K \colon (\m{C},P) \to (\m{E},R) \in \Mor(\FA_{\m{L}}) $ where $ K $ is full, surjective on objects and for each $ e \in \Ob(e) $, $ K^{p}_{e} $ is surjective, then there exists a unique morphism $ H \colon (\m{E},R) \to (\m{D},Q) $ such that $ H K = F $ if and only if $ \ker K \leq \ker F $. 
		\end{enumerate}
	\end{thrm}
	\begin{proof}[Proof of Part 1 of Theorem \ref{thm: Fibered homomorphism theorem}]
		Let $ T $ be the theory of $ (\m{C},P) $ and $ S $ the internal structure of $ (\m{C},P) $. Let $ T^{\prime} $ be the theory of $ F(S) $ (which is an $ \m{L}^{m} $-theory) and $ G $ the generic $ T^{\prime} $-model in $ (\m{C}_{T^{\prime}},P_{T^{\prime}}) $. Let $ \iota \colon (\m{C},P) \to (\m{C}_{T},P_{T}) $ be the morphism which with $ \overline{S} $ witnesses the equivalence of $ (\m{C},P) $ and $ (\m{C}_{T},P_{T}) $. We have the following diagram:
		\begin{equation*}
		\begin{tikzcd}
		(\m{C},P) \ar[swap]{d}{\textstyle \iota} \ar{rr}{\textstyle F} & & (\m{D},Q) \\
		(\m{C}_{T},P_{T}) \ar{r}{\textstyle \overline{G}} & (\m{C}_{T^{\prime}},P_{T^{\prime}}) \ar[swap]{ru}{ \textstyle \overline{F(S)}} & 
		\end{tikzcd}
		\end{equation*}
		It is straightforward to verify that the above diagram commutes. Since $ \overline{G} \circ \iota $ is injective on objects, from Proposition \ref{prop: The image of an injective on objects morphism determines a sub-prop-category.}, the following diagram commutes
		\begin{equation*}
		\begin{tikzcd}
		(\m{C},P) \ar[bend right,swap]{ddr}{\textstyle H} \ar{rr}{\textstyle F} \ar[swap]{dr}{\textstyle \overline{G}\circ \iota} & & (\m{D},Q) \\
		& (\m{C}_{T^{\prime}},P_{T^{\prime}})  \ar[swap]{ru}{\textstyle \overline{F(S)}} &  \\
		& (\overline{G} \circ \iota(\m{C}),\overline{G} \circ \iota(P)) \ar[hookrightarrow]{u}{\textstyle \lambda} \ar[swap, bend right]{ruu}{\textstyle \overline{F(S)} \circ \lambda} &
		\end{tikzcd}
		\end{equation*}
		where $ \lambda $ is the inclusion morphism and $ H $ the corestriction. Let $ \epsilon = H $, $ \psi = \overline{ F(S) } \circ \lambda$ and $ (\overline{G} \circ \iota(\m{C}), \overline{G} \circ \iota(P)) = (\m{E},R) $. By construction, $ \epsilon $ and $ \psi $ meet the conditions of the theorem. Since $ \iota $ is a sub-prop-morphism, from Proposition \ref{prop: FAL is closed under subprop-categories}, $ (\m{E},R) \in \FA_{\m{L}} $. 
	\end{proof}

	\begin{proof}[Proof of part 2 of Theorem \ref{thm: Fibered homomorphism theorem}]
		($ \impliedby $) Let $ e_{1} \in \Ob(\m{E}) $. Since $ K^{o} $ is surjective on objects, there exists $ c_{1} \in \Ob(\m{C}) $, such that $ K^{o}(c_{1}) = e_{1} $. Define $ H^{o}e_{1} \coloneqq F^{o}c_{1} $. Let $ e_{2} \in \Ob(\m{E}) $, and $ h_{1} \colon e_{1} \to e_{2} $. Then there exists $ c_{2} \in \Ob(\m{C}) $ such that $ K^{o}c_{2} = e_{2} $ and since $ K^{o} $ is full, there exists $ f_{1} \colon c_{1} \to c_{2} $ such that $ K^{o}f_{1} = h_{1} $. Define $ H^{o}h_{1} \coloneqq F^{o}f_{1} $. Since $ \ker K \leq \ker F $, $ H^{o} $ is well defined. If $ h_{2} \colon e_{2} \to e_{3} $, let $ f_{2} \colon c_{2} \to c_{3} $ such that $ K^{o}f_{2} = h_{2} $. Then 
		$$
		H^{o}(h_{2}h_{1}) = H^{o}(K^{o}f_{2}K^{o}f_{1}) = H^{o}(K^{o}(f_{2}f_{1})) = F^{o}(f_{2}f_{1}) = H^{o}h_{2}H^{o}h_{1}.
		$$
		Also, $ H^{o}(id_{e_{1}}) = H^{o}(K^{o}id_{c_{1}}) = F^{o}id_{c_{1}} = id_{F^{o}c_{1}}= id_{H^{o}e_{1}} $ and so $ H^{o} $ is a functor.
		
		Let $ 1_{\m{E}} $ be the terminal element in $ \m{E} $. Then there exists $ c \in \Ob(\m{C}) $ such that $ K^{o}c = 1_{\m{E}} $. Since $ K^{o} $ preserves finite products $ K^{o}1_{\m{C}} $ is a terminal element in $ \m{E} $ and so $ K^{o}1_{\m{C}} \cong 1_{\m{E}} $. Since $ H^{o} $ is a functor, it preserves isomorphisms and so $ H^{o}1_{\m{E}} \cong H^{o}K^{o}1_{\m{C}} = F^{o}1_{\m{C}}  $, and so $ H^{o}1_{\m{E}} $ is terminal in $ \m{D} $ since $ F^{o} $ preserves finite products.
		Let $ e_{1},e_{2} \in \Ob(\m{E}) $ and $ c_{1},c_{2},c \in \Ob(\m{C}) $, such that $ K^{o}c_{i} = e_{i} $ and $ K^{o}c = e_{1} \times e_{2} $. For $ i \in \{ 1,2 \} $, let $ f_{i} \colon c \to c_{i} $ such that $ K^{o}f_{i} = \pi_{i}^{e_{1},e_{2}} $ and let $ \alpha = \langle K^{o}\pi_{1}^{c_{1},c_{2}}, K^{o}\pi_{2}^{c_{1},c_{2}} \rangle $. Since $ K^{o} $ preserves finite products, $ \alpha $ is an isomorphism. For $ i \in \{ 1,2 \} $, $K^{o}\pi_{i}^{c_{1},c_{2}} K^{o}\langle f_{1},f_{2} \rangle = K^{o}f_{i} = \pi_{i}^{e_{1},e_{2}} $ and so $ K^{o}\langle f_{1},f_{2} \rangle = \alpha^{-1} $. Then $ H^{o}(\alpha^{-1}) = F^{o}\langle f_{1},f_{2} \rangle $ is an isomorphism and for $ i \in \{ 1,2 \} $, $ F^{o}\pi_{i}^{c_{1},c_{2}} \circ F^{o}\langle f_{1},f_{2} \rangle = F^{o}f_{i} $. Thus since $ F^{o} $ preserves finite products, $ F^{o}\langle f_{1}, f_{2} \rangle $ is a change-of-product isomorphism from $ (F^{o}c,F^{o}f_{1},F^{o}f_{2}) $ to $ (F^{o}(c_{1} \times c_{2}), F^{o}\pi_{1}^{c_{1},c_{2}},F^{o}\pi_{2}^{c_{1},c_{2}}) $. It follows that $ (H^{o}(e_{1} \times e_{2}), H^{o}\pi_{1}^{e_{1},e_{2}},H^{o}\pi_{2}^{e_{1},e_{2}})= (F^{o}c,F^{o}f_{1},F^{o}f_{2}) $ is a product diagram for $ H^{o}e_{1}, H^{o}e_{2} $ and so $ H^{o} $ preserves finite products.
		
		Let $ e \in \Ob(\m{E}) $, and $ r \in R(e) $. Then there exists $ c_{1} \in \Ob(\m{C}) $ and $ r_{1} \in P(c_{1}) $ such that $ K^{p}_{c_{1}}(r_{1}) = r $. Define $ H^{p}_{e}(r) \coloneqq F^{p}_{c_{1}}(r_{1}) $. If $ r_{2} \in P(c_{2}) $, such that $ K^{p}_{c_{2}}(r_{2}) = r $, then $ (r_{1},r_{2}),(r_{2},r_{1}) \in \ker K \subseteq \ker F $ and so $ F^{p}_{c_{1}}(r_{1}) = F^{p}_{c_{2}}(r_{2}) $. It follows that $ H^{p}_{e} $ is well-defined and by definition, $ H^{p}_{e} \circ K^{p}_{c_{1}} = F^{p}_{c_{1}} $. Let $ r_{1},r_{2} \in R(e) $ such that $ r_{1} \leq r_{2} $. Then there exists $ c \in \Ob(\m{C}) $ and $ r_{1}^{\prime},r_{2}^{\prime} \in P(c) $ such that $ K^{p}_{c}(r_{i}^{\prime}) = r_{i} $ for $ i \in \{ 1,2 \} $. Then $ (r_{1}^{\prime},r_{2}^{\prime}) \in \ker K \subseteq \ker F $, and so $ H^{p}_{e}(r_{1}) = F^{p}_{c}(r_{1}^{\prime}) \leq F^{p}_{c}(r_{2}^{\prime}) = H^{p}_{e}(r_{2}) $. Thus $ H^{p}_{e} $ is monotone.

		Let $ h \colon e_{1} \to e_{2} $ and $ f \colon c_{1} \to c_{2} $ such that $ K^{o}f = h $. Then
		\begin{align*}
		Q(H^{o}h) \circ H^{p}_{e_{2}} \circ K_{c_{2}}^{p} & = Q(F^{o}f) \circ F^{p}_{c_{2}} = F^{p}_{c_{1}} \circ P(f) \\
		& = H_{e_{1}}^{p} \circ K_{c_{1}}^{p} \circ P(f) = H_{e_{1}}^{p} \circ R(h) \circ K_{c_{2}}^{p},
		\end{align*}
		and since $ K^{p}_{c_{2}} $ is surjective, $ Q(H^{o}h) \circ H^{p}_{e_{2}} = H^{p}_{e_{1}} \circ R(h) $. Thus, $ H^{p} \colon R \Rightarrow QH^{o} $ is a natural transformation.

		Let $ e \in \Ob(\m{E}) $. Then there exists $ c,c^{\prime} \in \Ob(\m{C}) $, such that $ K^{o}c^{\prime} = e \times e $ and $ K^{o}c = e $. Then for $ i \in \{ 1,2 \} $, there exists $ f_{i} \colon c^{\prime} \to c $ such that $ K^{o}f_{i} = \pi_{i}^{e,e} $. Then for $ i \in \{ 1,2 \} $, 
		$$ 
		K^{o}\pi_{i}^{c,c} K^{o} \langle f_{1},f_{2} \rangle = K^{o}f_{i} = \pi_{i}^{e,e} .
		$$
		It follow that $ K^{o}\langle f_{1},f_{2} \rangle = \alpha^{-1} $, where $ \alpha = \langle K^{o}\pi_{1}^{c,c},K^{o}\pi_{2}^{c,c} \rangle $. Define:
		\begin{equation*}
		a \coloneqq \langle F^{o}\pi_{1}^{c,c},F^{o}\pi_{2}^{c,c} \rangle \quad b \coloneqq \langle H^{o}\pi^{e,e}_{1}, H^{o}\pi_{2}^{e,e}\rangle.
		\end{equation*}
		Then
		\begin{equation*}
		a F^{o}\langle f_{1},f_{2} \rangle = \langle F^{o}f_{1},F^{o}f_{2} \rangle = \langle H^{o} \pi_{1}^{e,e} ,H^{o}\pi_{2}^{e,e} \rangle = b,
		\end{equation*}
		and so $ F^{o}\langle f_{1},f_{2} \rangle = a^{-1}b $.
		Also note that
		\begin{equation*}
		F^{p}_{c \times c}(Eq_{c}) = Q(a)(Eq_{F^{o}c}) \quad \ \text{and} \quad K_{c \times c}^{p}(Eq_{c}) = R(\alpha)(Eq_{K^{o}c}).
		\end{equation*}
		Then
		\begin{align*}
		Q(b)(Eq_{H^{o}e}) & = Q(b)(Eq_{F^{o}c}) = Q(b)(Q(a^{-1})F_{c \times c}^{p}Eq_{c}) \\
		& = Q(a^{-1}b)F_{c \times c}^{p}(Eq_{c}) = Q(F^{o} \langle f_{1},f_{2}\rangle)F^{p}_{c \times c}(Eq_{c}) \\
		& = F^{p}_{c^{\prime}}  P(\langle f_{1},f_{2} \rangle)(Eq_{c}) = H^{p}_{e \times e}K^{p}_{c^{\prime}}P(\langle f_{1},f_{2} \rangle)(Eq_{c}) \\
		& = H^{p}_{e \times e}R(K^{o}\langle f_{1},f_{2}\rangle)K^{p}_{c \times c}(Eq_{c}) = H_{e \times e}^{p}(R(\alpha^{-1})K^{p}_{c \times c}(Eq_{c})) \\
		& = H^{p}_{e \times e}(Eq_{K^{o}c}) = H^{p}_{e \times e}(Eq_{e}).
		\end{align*}
		
		Let $ e_{1},e_{2} \in \Ob(\m{E}) $, and $ c_{1},c_{2},c \in \Ob(\m{C}) $ such that $ K^{o}c_{i} = e_{i} $ and $ K^{o}c = e_{1} \times e_{2} $. Let $ f_{i} \colon c \to c_{i} $ such that $ K^{o}f_{i} = \pi_{i}^{e_{1},e_{2}} $ and define
		\begin{equation*}
		a \coloneqq \langle F^{o}\pi_{1}^{c_{1},c_{2}},F^{o}\pi_{2}^{c_{1},c_{2}} \rangle \quad b \coloneqq \langle H^{o}\pi^{e_{1},e_{2}}_{1}, H^{o}\pi_{2}^{e_{1},e_{2}}\rangle.
		\end{equation*}
		Also, let $ \alpha \coloneqq \langle K^{o}\pi_{1}^{c_{1},c_{2}},K^{o}\pi_{2}^{c_{1},c_{2}}\rangle $. Then for $ \Omega \in \s{L}_{q} $
		\begin{align*}
		& = \Omega_{H^{o}e_{1},H^{o}e_{2}}\circ Q(b^{-1}) \circ H_{e_{1} \times e_{2}}^{p} \circ R(\alpha^{-1}) \circ K_{c_{1} \times c_{2}}^{p} \\
		& = \Omega_{F^{o}c_{1},F^{o}c_{2}} \circ Q(a^{-1})\circ Q(b^{-1}a) \circ H_{e_{1} \times e_{2}}^{p} \circ R(\alpha^{-1}) \circ K_{c_{1} \times c_{2}}^{p} \\
		& = \Omega_{F^{o}c_{1},F^{o}c_{2}} \circ Q(a^{-1})\circ Q(H^{o}\alpha) \circ H_{e_{1} \times e_{2}}^{p} \circ R(\alpha^{-1}) \circ K_{c_{1} \times c_{2}}^{p} \\
		& = \Omega_{F^{o}c_{1},F^{o}c_{2}} \circ Q(a^{-1})\circ H^{p}_{K^{o}(c_{1} \times c_{2})} \circ R(\alpha) \circ R(\alpha^{-1}) \circ K_{c_{1} \times c_{2}}^{p} \\
		& = \Omega_{F^{o}c_{1},F^{o}c_{2}} \circ Q(a^{-1})\circ H^{p}_{K^{o}(c_{1} \times c_{2})} \circ K_{c_{1} \times c_{2}}^{p} \\
		& = \Omega_{F^{o}c_{1},F^{o}c_{2}} \circ Q(a^{-1}) \circ F^{p}_{c_{1} \times c_{2}} \\
		& = F^{p}_{c_{1}} \circ \Omega_{c_{1},c_{2}} \\
		& = H^{p}_{e_{1}} \circ K^{p}_{c_{1}} \circ \Omega_{c_{1},c_{2}} \\
		& = H^{p}_{e_{1}} \circ \Omega_{K^{o}c_{1},K^{o}c_{2}} \circ R(\alpha^{-1}) \circ K^{p}_{c_{1} \times c_{2}} \\
		& = H^{p}_{e_{1}} \circ \Omega_{e_{1},e_{2}} \circ R(\alpha^{-1}) \circ K^{p}_{c_{1} \times c_{2}}.
		\end{align*}
		Since $ R(\alpha^{-1}) \circ K^{p}_{c_{1} \times c_{2}} $ is surjective it follows that
		$$
		\Omega_{H^{o}e_{1},H^{o}e_{2}}\circ Q(b^{-1}) \circ H_{e_{1} \times e_{2}}^{p} = H^{p}_{e_{1}} \circ \Omega_{e_{1},e_{2}}.
		$$
		Thus $ H $ is a morphism in $ \FA $ and $ H\circ K = F $.
		
		($ \implies $). Since $ \m{K} $ is full, surjective on objects and for each $ e \in \Ob(\m{E}) $, $ K^{p}_{e} $ is surjective, there is only one possible definition for $ H $. One may verify that if $ \ker K \nleq \ker F $, then either $ H^{o} $ is ill-defined or for some $ e \in \Ob(\m{E}) $, $ H^{p}_{e} $ is ill-defined or not monotone.
	\end{proof}
	
	\begin{rem}
		Like part 1 of Theorem \ref{thm: Fibered homomorphism theorem}, the internal logic can also be used to prove part 2. In what follows we sketch the proof of the ( $ \impliedby $) direction.
		Let $ S $ be the internal structure in $ (\m{E},R) $ and $ T $ its theory. We construct a $ T $-model $ A \in (\m{D},Q) $, so that
		$$
		\begin{tikzcd}
		(\m{C},P) \ar[swap]{dr}{\textstyle K} \ar{rr}{\textstyle F} & & ( \m{D},Q ) & \\
		& (\m{E},R) \ar[swap,dashed]{ru}{\textstyle \overline{A} \circ \iota } \ar[swap]{rr}{\textstyle \iota} & & (\m{C}_{T},P_{T}) \ar[swap]{ul}{\textstyle \overline{A}}
		\end{tikzcd}
		$$
		commutes.
		
		For each $ e \in \Ob(\m{E}) $, let $ c_{e} \in \Ob(\m{C}) $ such that $ K^{o}c_{e} = e $.  For each $ h \in \Mor(\m{E}) $, let $ f_{h} \in \Mor(\m{C}) $, such that $ K^{o}f_{h} = h $ and for each $ e \in \Ob(\m{E}) $ and $ r \in R(e) $, let $ q_{r} \in P(c) $ for some $ c \in \Ob(\m{C}) $ such that $ K^{p}_{c}(q_{r}) = r $. For each each sort $ e $, $ A\db{e} \coloneqq F^{o}c_{e} $, for each unary function symbol $ h \colon e_{1} \times \cdots \times e_{n} \to e $, $ A \db{h} \coloneqq F^{o}f_{h} $ and each unary relation symbol $ r \in R(e_{1} \times \cdots \times e_{n}) $, $ A \db{r} \coloneqq F^{p}_{c_{\bar{e}}}q_{r} $. For each $ n $-ary function symbol $ h \colon e_{1}, \ldots, e_{n} \to e $, we define $ A \db{h} \coloneqq f_{h} \circ a_{\bar{e}}^{-1} $ where $ a_{\bar{e}} \colon F^{o}(c_{\bar{e}}) \to A\db{e_{1}} \times \cdots \times A \db{e_{n}} $ is the change-in-product isomorphism. Similarly, for each $ n $-ary relation symbol $ r \subseteq e_{1}, \ldots, e_{n} $, we define $ A \db{r} \coloneqq Q(a_{\bar{e}})^{-1}\circ F^{p}_{c_{\bar{e}}}(q_{r}) $.
		
		By induction, one may prove that for each term $ M : e \ [\Gamma] $, where $ \Gamma = x_{1}:e_{1}, \ldots, x_{n}:e_{n} $, that $ A \db{M: e \ [\Gamma]} = F^{o}f_{S \db{M: e \, [\Gamma]}} \circ a_{\bar{e}}^{-1} $, and for each formula-in-context $ \phi \ [\Gamma] $, that $ A \db{\phi [\Gamma]}= Q(a^{-1}_{\bar{e}})\circ F^{p}_{c_{\bar{e}}}(q_{S \db{\phi [\Gamma]}}) $. It follows that $ A $ satisfies each equation-in-context in $ T $, and if 	$ \phi_{1}, \ldots, \phi_{n} \vdash \phi_{n+1} \ [\Gamma] $ is in $ T $, then 
		\begin{align*}
		& S \db{\phi_{1} \con \ldots \con \phi_{n} \ [\Gamma]} \leq S \db{\phi_{n+1}[\Gamma]} \\
		\implies \ & K^{p}_{c}(q_{S \db{\phi_{1} \con \ldots \con \phi_{n}[\Gamma]}}) \leq K^{p}_{c}(q_{S \db{\phi_{n+1}[\Gamma]}})  \\
		\implies \ & F^{p}_{c}(q_{S \db{\phi_{1} \con \ldots \con \phi_{n}[\Gamma]}}) \leq F^{p}_{c}(q_{S \db{\phi_{n+1}[\Gamma]}}) \\
		\implies \ & F^{p}_{c}(q_{S\db{\phi_{1}[\Gamma]}}) \con \ldots \con F_{c}^{p}(q_{S \db{\phi_{n}[\Gamma]}}) \leq F^{p}_{c}(q_{S \db{\phi_{n+1}[\Gamma]}}) \\
		\implies \ & A \db{\phi_{1}[\Gamma]} \con \ldots \con A \db{\phi_{n}[\Gamma]} \leq A \db{\phi_{n+1}[\Gamma]}.
		\end{align*}
		It follows that $ A $ is a $ T $-model. By construction, we have $ K \circ \overline{A} \circ \iota = F $.
	\end{rem}
	
	\begin{prop}
		The classes $ \s{E} $ and $ \s{M} $ form a factorization system in $ \FA_{\m{L}} $. 
	\end{prop}
	\begin{proof}
		Note that both $ \s{M} $ and $ \m{E} $ contain all the isomorphisms in $ \FA_{\m{L}} $ and are closed under composition. From Part 1 of Theorem \ref{thm: Fibered homomorphism theorem}, each $ F \in \Mor(\FA_{\m{L}}) $ can be factored as $ m \circ e $ with $ m \in \s{M} $ and $ e \in \s{E} $. Now consider the following commuting solid diagram of morphisms in $ \FA_{\m{L}} $:
		\begin{equation*}
		\begin{tikzcd}
		\cdot \ar[swap]{d}{e} \ar{r}{u} & \cdot \ar{d}{e^{\prime}} \\
		\cdot \ar[swap]{d}{m} \ar[dashed]{r}{\exists ! w} & \cdot \ar{d}{m^{\prime}} \\
		\cdot \ar{r}{v} & \cdot
		\end{tikzcd}
		\end{equation*}
	We want to to show there is a unique morphism $ w $ making the small squares commute. Note that $ \Th(e) \subseteq \Th(vme) = \Th(m^{\prime}e^{\prime}u) = \Th(e^{\prime}u) $, and since $ e $ is bijective on objects, $ \ker(e) \leq \ker(e^{\prime}u) $. From part 2 of Theorem \ref{thm: Fibered homomorphism theorem}, there exists a unique morphism $ w $ such that $ w e = e^{\prime} u $. Moreover it is straightforward to verify that the morphisms in $ \s{E} $ are epimorphisms in $ \FA_{\m{L}} $ and so $ e $ is an epimorphism. From this fact and the fact that the large rectangle and small top square commute, the small bottom square must commute as well.
	\end{proof}

	The above factorization system factors morphisms into their ``logical'' ($ \s{E} $) and ``nonlogical'' $ (\s{M}) $ parts. Given a morphism $ F \colon (\m{C},P) \to (\m{D},Q) $, We can define the ``theory of $ F $'', $ \Th(F) $, to be $ \Th(F \circ \overline{S}) $, where $ \overline{S} \colon (\m{C}_{T},P_{T}) \to (\m{C},P) $ is the canonical morphism induced by the internal structure $ S $ of $ (\m{C},P) $. If $ m \circ \epsilon $ is an $ \s{E}\s{M} $-factorization of $ F $, then $ \Th(F) = \Th(\epsilon) $, and $ \epsilon $ is universal in the sense that for each $ H \colon (\m{C},P) \to (\m{E},R) $ such that $ \Th(H) \geq \Th(F) $, there exist a unique morphism $ K $ such that $ H = K \circ \epsilon $.

	\section{Algebraic Characterizations of Logical Closure Operators}\label{sec: Semantic Closure Operators}
	
	Let $ Sg $ be a single-sorted algebraic signature, $ \Alg_{Sg} $ the class of set-valued $ Sg $-algebras and $ Eq_{Sg} $ the class of $ Sg $-equations. Then there is an adjunction
	\begin{equation}\label{eqn: Classical Birkhoff Adjunction}
	\begin{tikzcd}
	\s{P}(Eq_{Sg})
	\arrow[r, "\Alg(\cdot)"{name=F}, bend left=25] &
	\s{P}(\Alg_{Sg})^{op}
	\arrow[l, "Eq(\cdot)"{name=G}, bend left=25]
	\arrow[phantom, from=F, to=G, "\dashv" rotate=-90]
	\end{tikzcd}
	\end{equation}
	Where $ \s{P} $ is the operation of taking the class of all subclasses and both $ \s{P}(Eq_{Sg}) $ and $ \s{P}(\Alg_{Sg}) $ are ordered by inclusion. $ \Alg(\cdot) $ takes a class of algebras to the collection of all $ Sg $-equations they satisfy and $ Eq(\cdot) $ takes a collection of $ Sg $-equations to the class of all $ Sg $-algebras that satisfy them. This adjunction determines a closure operator $ Eq \circ \Alg $ on $ \s{P}(Eq_{Sg}) $, and a closure operator $ \Alg \circ Eq $ on $ \s{P}(\Alg_{Sg}) $. From Birkhoff's Completeness Theorem for Equational Logic, $ Eq \circ \Alg(\Theta) $ is the closure of $ \Theta $ under the derivation rules of equational logic. Birkhoff's HSP Theorem asserts $ \Alg \circ Eq(\cdot) $ is $ \mathbb{H}\mathbb{S}\mathbb{P}(\cdot) $  where $ \mathbb{H} $, $ \mathbb{S} $, $ \mathbb{P} $ are the operations which close a class of algebras under the operations of taking homomorphic images, subalgebras and products respectively \cite{Birkhoff1935}.
	
	In the categorical semantics for equational logic, $ \Alg_{Sg}$ is equivalent to $ \FP(\m{C}_{Sg},\Set) $ where $ \FP $ is the $ 2 $-category of categories with finite products and product preserving functors. Birkhoff's HSP Theorem interpreted in the categorical semantics says $ \ca{Y} \subseteq \Ob(\FP(\m{C}_{Sg},\Set)) $ is an equational class iff $ \mathbb{H}\mathbb{S}\mathbb{P}(\m{Y}) = \m{Y} $, where
	\begin{enumerate}
		\item $ \mathbb{P}(\m{Y}) $ is the smallest class containing $ \m{Y} $ stable under products.
		\item $ \mathbb{S}(\m{Y}) $ is the smallest class containing $ \m{Y} $ such that whenever $ F \in \m{Y} $ and $ \eta \colon K \Rightarrow F $ is a monomorphism, then $ K \in \mathbb{S}(\m{Y}) $.
		\item $ \mathbb{H}(\m{Y}) $ is the smallest class containing $ \m{Y} $ such that whenever $ F \in \m{Y} $, and $ \eta \colon F \Rightarrow K $ is a regular epimorphism, then $ K \in \mathbb{H}(\m{Y}) $.
	\end{enumerate}
	For an analogous characterization for a multi-sorted algebraic signature see \cite{Adamek2011}. Since we are considering first-order logics, we are interested in a similar result where $ \s{P}(Eq_{Sg}) $ is replaced by the complete lattice of all $ Sg $-theories $ \Th_{Sg} $ and algebras by first-order structures. If $ Sg $ is single-sorted and $\s{L}$ is a language for classical first-order logic, then the Tarskian $Sg$-structures of classical first-order logic can be identified with morphisms $ F \colon (\m{C}_{Sg},P_{Sg}) \to (\Set, \s{P}) $, where $ \s{P} $ is the preimage functor. It is well known that closing a class of Tarskian structures under their common theory, is equivalent to closing the class under ultraproducts, isomorphic copies and ultraroots \cite[p.~454]{Hodges1993}. In what follows, we develop an analogous result for an arbitrary first-order logic $ \m{L} $ and its general prop-categorical semantics. 
	
	Let $ \Mod_{Sg}^{\m{L}} $ be the collection of all $ Sg $-structures in $ \FA_{\m{L}} $. Then we are interested in characterizing  $ \Mod_{(\cdot)} \circ \Th $, where
	$$
	\begin{tikzcd}
	{\Th_{Sg}}
	\arrow[r, pos = .6, "\Mod_{(\cdot)}"{name=F}, bend left=25] &
	\s{P}(\Mod_{Sg}^{\m{L}})^{op}.
	\arrow[l, pos = .4, "\Th(\cdot)"{name=G}, bend left=25]
	\arrow[phantom, from=F, to=G, "\dashv" rotate=-90]
	\end{tikzcd}
	$$
	As before, we identify $ \Mod_{Sg}^{\m{L}} $ with the collection of all morphisms $ F \colon (\m{C}_{Sg},P_{Sg}) \to (\m{C},P) $ such that $ (\m{C},P) \in \Ob(\FA_{\m{L}}) $. Then, for each $ \m{Y} \subseteq \Mod_{Sg}^{\m{L}} $, $ \Th(\m{Y}) $ is the $ Sg $-theory whose assertions are $ \{ a \in A_{Sg} : \forall F \in \m{Y}, G\db{a} \in \ker F \} $ and for each $ T \in \Th_{Sg}$, $ \Mod_{T} = \{ F \in \Mod_{Sg}^{\m{L}}: G \db{T} \subseteq \ker F \}  $. If $ \FA_{\m{L}} $ forms a complete semantics for $ \m{L} $, $ \Th( \Mod_{T})$ is the $Sg $-theory whose assertions are $  \{ a \in \A_{Sg} : T \vdash_{\m{L}} a \} $.
	For $ \m{Y} \subseteq \Mod_{Sg}^{\m{L}} $, we wish to characterize $ \Mod_{(\cdot)} \circ \Th(\m{Y}) $. Let $ \{ F_{i} \colon (\m{C}_{Sg},P_{Sg}) \to (\m{C}_{i},P_{i})  \}_{i \in I} \subseteq \Mod_{(\cdot)} \circ \Th(\m{Y}) $. Then each $ (\m{C}_{i},P_{i}) \in \Ob(\FA_{\m{L}}) $, and from Proposition \ref{prop: Logic of the product is stronger than the logic of the set.}, $ (\prod \m{C}_{i},\prod P_{i}) \in \Ob(\FA_{\m{L}}) $, and so $ \langle F_{i} \rangle_{i \in I}  \in \Mod_{Sg}^{\m{L}} $. Moreover, $ \langle F_{i} \rangle_{i \in I} $ satisfies each assertion satisfied by all $ F_{i} $ and so $ \langle F_{i} \rangle_{i \in I} \in \Mod_{(\cdot)} \circ \Th(\m{Y}) $. We define the \textbf{(external) product of $ \{ F_{i} \}_{i \in I} $} to be $ \langle F_{i} \rangle_{i \in I} $ and let $ \mathbb{P}(\m{Y}) $ denote the closure of $ \m{Y} $ under external products. And so, $ \Mod_{(\cdot)} \circ \Th(\m{Y}) $ is stable under taking external products.
	
	Let $ F \colon (\m{C}_{Sg},P_{Sg}) \to (\m{C},P) \in \Mod_{(\cdot)} \circ \Th(\m{Y}) $. For each $ (\m{D},Q) \in \Ob(\FA_{\m{L}}) $ and $ H \in \FA_{\m{L}}((\m{C},P),(\m{D},Q)) $, $ H \circ F \in \Mod_{(\cdot)} \circ \Th(\m{Y}) $. We say $ H \circ F $ is an \textbf{(external) homomorphic image of $ F $} and let $ \mathbb{H}(\m{Y}) $ denote the closure of $ \m{Y} $ under external homomorphic images. If $ H \colon (\m{C}_{Sg},P_{Sg}) \to (\m{E},R) $ and $ \iota \colon (\m{E},R) \to (\m{C},P) $ are morphisms in $ \FA_{\m{L}} $ such that $ F = \iota \circ H $ and $ \iota $ is a sub-prop-morphism, then we call $ H $ an \textbf{(external) submodel of $ F $} and let $ \mathbb{S}(\m{Y}) $ denote the closure of $ \m{Y} $ under taking external submodels. Since $ G\db{\Th(\m{Y})} \subseteq \ker F $, $ \iota \circ H = F $ and $ \iota $ is a sub-prop-morphism, $ G \db{\Th(\m{Y})} \subseteq \ker H $ and so $ H \in \Mod_{(\cdot)} \circ \Th(\m{Y}) $. 
	
	\begin{rem}\label{rem: Comparing internal and external products}
		If we consider algebras satisfying some equational theory $ T $ as product preserving functors $ \{ F_{i} \colon \m{C}_{T} \to \m{C} \}_{i \in I} $, where $ \m{C} $ has arbitrary products then,
		\begin{equation*}
		\adjustbox{scale=1.1,center}{
			\begin{tikzcd}
			\m{C}_{T} \arrow[rr, "F_{\times}"{name = {P}}] \ar[swap]{dr}{\langle F_{i}\rangle_{i \in I}} & & \m{C}	\\
			& \prod_{i \in I}\m{C}  \ar[swap]{ru}{\times}  \arrow[Rightarrow,from = P, shorten = 5, swap, "a"] & 
			\end{tikzcd}}
		\end{equation*}
		commutes up to a change-in-product natural isomorphism $ a $, where $ \times $ is the right adjoint to the diagonal functor $ \Diag \colon \m{C} \to \prod_{i \in I} \m{C} $ and $ F_{\times} $ is the usual ``internal'' product of the algebras $ \{ F_{i} \}_{i \in I} $. Since $ \times $ is faithful, $ \langle F_{i} \rangle_{i \in I} $ and $ F_{\times} $ satisfy the same equations. For first-order models, the internal product need not correspond to the external. The issue is that for $ c \in \Ob(\prod_{i \in I}\m{C}) $, in general, $ \prod_{i \in I}P(c_{i}) \not\cong P(\times(c))  $.
	\end{rem}
	
	\begin{thrm}\label{thm: Characterization of Theory Closure.}
		Let $ \m{L} $ be a logic and $ Sg $ a small signature (not a proper class). For each $ \m{Y} \subseteq \Mod_{Sg}^{\m{L}} $, $ \Mod_{(\cdot)} \circ \Th(\m{Y}) = \mathbb{H} \mathbb{S}\mathbb{P}(\m{Y}) $.
	\end{thrm}
	\begin{proof}
		If $ \m{Y} $ is the class of all $ Sg $-$ \m{L} $ structures of a given theory $ T $, then we showed $ \m{Y} $ is closed under taking (external) products, submodels and homomorphic images and so $ \mathbb{H}\mathbb{S}\mathbb{P}(\m{Y}) \subseteq \Mod_{(\cdot)} \circ \Th(\m{Y}) $.
		
		In the other direction, let $ T = \Th(\m{Y}) $, $ F \colon (\m{C}_{Sg},P_{Sg}) \to (\m{C},P) $ be a $ T $-model, where $ (\m{C},P) \in \Ob(\FA_{\m{L}}) $ and let $ \{ T_{i} \}_{i \in I} $ be the collection of all $ Sg $-$ \m{L} $ theories such that there exists $ F_{i} \colon (\m{C}_{Sg},P_{Sg}) \to (\m{C}_{i},P_{i}) \in \m{Y} $ such that $ \Th(F_{i}) = T_{i} $. Then $ T = \bigcap_{i \in I}T_{i} $, and define $ \epsilon_{T} \coloneqq \overline{G}  \colon (\m{C}_{Sg},P_{Sg}) \to (\m{C}_T,P_{T}) $, where $ G \in (\m{C}_{T},P_{T}) $ is the generic $ T $-model. Then $ \ker \epsilon_{T} \leq \ker \langle F_{i} \rangle_{I}  $, and so from Theorem \ref{thm: Fibered homomorphism theorem}, there exists $ \iota \colon (\m{C}_{T},P_{T}) \to (\prod_{I} \m{C}_{i},\prod_{I} P_{i}) $ such that 
		\begin{equation*}
		\begin{tikzcd}
		(\m{C}_{Sg},P_{Sg}) \ar[swap]{dr}{\textstyle \epsilon_{T}} \ar{rr}{\textstyle \langle F_{i}\rangle_{I}} & & (\prod_{I} \m{C}_{i}, \prod_{I} P_{i})\\
		& (\m{C}_{T},P_{T}) \ar[swap,dashed]{ru}{\textstyle \iota } &
		\end{tikzcd}
		\end{equation*}
		commutes. Moreover, since $ \Th (\langle F_{i} \rangle_{I}) = T  $, $ \iota $ is a sub-prop-morphism. From Proposition \ref{prop: FAL is closed under subprop-categories}, $ (\m{C}_{T},P_{T}) \in \Ob(\FA_{\m{L}}) $ and so $ \epsilon_{T} \in \mathbb{S}\mathbb{P}(\m{Y}) = \m{Y} $. Since $ F $ is a $ T $-model, $ \ker \epsilon_{T} \leq \ker F $.  From Theorem \ref{thm: Fibered homomorphism theorem}, there exists a morphism $ H \colon (\m{C}_{T},P_{T}) \to (\m{C},P) $, such that 
		\begin{equation*}
		\begin{tikzcd}
		(\m{C}_{Sg},P_{Sg}) \ar[swap]{dr}{\textstyle \epsilon_{T}} \ar{rr}{\textstyle F} & & ( \m{C},P)\\
		& (\m{C}_{T},P_{T}) \ar[swap,dashed]{ru}{\textstyle H } &
		\end{tikzcd}
		\end{equation*}
		commutes and so $ F \in \mathbb{H}\mathbb{S}\mathbb{P}(\m{Y}).$ Therefore, $ \mathbb{H} \mathbb{S} \mathbb{P}(\m{Y}) \supseteq \Mod_{(\cdot)} \circ \Th(\m{Y}) $.
	\end{proof}

	Let $ \Log $ be the partial order of all logics. We now consider the following adjunction:
	\begin{equation}
	\begin{tikzcd}
	\phantom{T} \Log
	\arrow[r, "\Ob(\FA_{(\cdot)})"{name=F}, bend left=25] &
	\s{P}(\Ob(\FA))^{op}
	\arrow[l, "\vDash_{(\cdot)}"{name=G}, bend left=25]
	\arrow[phantom, from=F, to=G, "\dashv" rotate=-90]
	\end{tikzcd}\label{eqn: signature independent adjunction}
	\end{equation} In order to give a characterization of the closure operator $ \Ob(\FA_{\vDash_{(\cdot)}}) $, we first restrict the logics in $ \Log$ to some fixed signature $ Sg $. Let $ \Log^{Sg} $ be the collection of all logics restricted to $ Sg $-assertions and we denote the corresponding restrictions of $ \Ob(\FA_{(\cdot)}) $ and $ \vDash_{(\cdot)} $, $ \Ob(\FA_{(\cdot)}^{Sg}) $ and $ \vDash_{(\cdot)}^{Sg} $ respectively:
	\begin{equation}\label{eqn: restricted adjunction.}
	\begin{tikzcd}
	\Log^{Sg}
	\arrow[r, "\Ob(\FA_{(\cdot)}^{Sg})"{name=F}, bend left=25] &
	\s{P}(\Ob(\FA)^{Sg})^{op}
	\arrow[l, "\vDash_{(\cdot)}^{Sg}"{name=G}, bend left=25]
	\arrow[phantom, from=F, to=G, "\dashv" rotate=-90]
	\end{tikzcd}
	\end{equation}
	
	In the context of (untyped) equational logic Adjunction \ref{eqn: restricted adjunction.} corresponds to the following: Let $ \s{L}_{\omega} $ be an algebraic signature and $ V $ a set of variables of cardinality $ \lambda $. Let $ Eq_{V} $ be the collection of $ \s{L}_\omega $-equations over $ V $ which we identify with $ Fm_{V}^{2} $, the square of the formula-algebra. A collection of $ \s{L}_{\omega} $-algebras $ \m{A} $ determines a structural closure operator $ \vDash_{\m{A}} $, defined by $ \Theta \vDash_{\m{A}} \epsilon = \delta $ if and only if for all homomorphisms $ h \colon Fm_{V} \to A $, where $ A \in \m{A} $, if $ \Theta \subseteq \ker h $, then $ \epsilon = \delta \in \ker h $. Observe that $ \vDash_{\mathbb{S}\mathbb{P}(\m{A})} \ = \ \vDash_{\ca{A}} $, but $ \mathbb{S}\mathbb{P}(\m{A}) $ may not be the largest collection of algebras defining the same consequence. In \cite{Blok2006} it is shown that the largest such class of $ \s{L}_{\omega} $-algebras is $ \mathbb{U}_{\lambda}\mathbb{S} \mathbb{P}(\m{A}) $, where $ B \in \mathbb{U}_{\lambda}(\m{A}) $ if  every $ \lambda $-generated subalgebra of $ B $ is in $ \m{A} $.
	
	Taking the operation $ \mathbb{U}_{\lambda} $ as inspiration, where $ \lambda $ up to renaming specifies the propositional signature, given a signature $ Sg $, and $ \m{X} $ a collection of prop-categories, we define $ \mathbb{U}_{Sg}(\ca{X}) $ so that $ (\m{D},Q) \in \mathbb{U}_{Sg}(\m{X}) $, if each classifying sub-prop-category $ (\m{C}_{T},P_{T})  $ of $ (\m{D},Q)  $ is in $ \m{X} $, where $ \Sg(T) = Sg $. For $ \m{X} \subseteq \Ob(\FA) $, we take $ \mathbb{P}(\m{X}) $ to be the closure of $ \m{X} $ under taking products and $ \mathbb{S}(\m{X}) $ to be the closure of $ \m{X} $ under taking subprop-categories. 
	
	\begin{thrm}\label{thm: Fixed signature adjunction}
		For $ \m{X} \subseteq \Ob(\FA) $, and $ Sg $ a small signature (not a proper class), $ \Ob(\FA^{Sg}_{(\cdot)}) \circ \vDash^{Sg}_{(\cdot)}(\m{X}) = \mathbb{U}_{Sg} \mathbb{S}\mathbb{P}(\m{X}) $. 
	\end{thrm}
	\begin{proof}
		From Proposition \ref{prop: FAL is closed under subprop-categories} and Proposition \ref{prop: Logic of the product is stronger than the logic of the set.}, $\Ob(\FA^{Sg}_{(\cdot)}) \circ \vDash^{Sg}_{(\cdot)}(\m{X}) $ is stable under $ \mathbb{S} $ and $ \mathbb{P} $. Let $ (\m{D},Q) \in \mathbb{U}_{Sg}(\Ob(\FA^{Sg}_{(\cdot)}) \circ \vDash^{Sg}_{(\cdot)}(\m{X})) $ and suppose $ T \vDash_{\m{X}} a $. Let $ F \colon (\m{C}_{Sg},P_{Sg}) \to (\m{D},Q) $ be a $ T $-model in $ (\m{D},Q) $ and $ T^{\prime} = \Th(F) $. From Theorem \ref{thm: Fibered homomorphism theorem} there exists $ \iota \colon (\m{C}_{T^{\prime}},P_{T^{\prime}}) \to (\m{D},Q) $ which makes the diagram
		\begin{equation*}
		\begin{tikzcd}
		(\m{C}_{Sg},P_{Sg}) \ar[swap]{dr}{\textstyle \epsilon_{T^{\prime}}} \ar{rr}{\textstyle F} & & ( \m{D},Q)\\
		& (\m{C}_{T^{\prime}},P_{T^{\prime}}) \ar[swap,dashed]{ru}{\textstyle \iota } &
		\end{tikzcd}
		\end{equation*}
		commute. Moreover, since $ \Th(F)= \Th(\epsilon_{T^{\prime}}) $, $ \iota $ is a sub-prop-morphism. Then $ (\m{C}_{T^{\prime}},P_{T^{\prime}}) \in \Ob(\FA^{Sg}_{(\cdot)}) \circ \vDash^{Sg}_{(\cdot)}(\m{X}) $, because $ (\m{D},Q) \in \mathbb{U}_{Sg}(\Ob(\FA^{Sg}_{(\cdot)}) \circ \vDash^{Sg}_{(\cdot)}(\m{X})) $. Since $ T \vDash_{\m{X}} a $ and $ G \db{T} \subseteq \ker \epsilon_{T^{\prime}} $, $ G \db{a} \in \ker \epsilon_{T^{\prime}} $. Thus $ G \db{a} \in \ker F $ and so $ T \vDash_{(\m{D},Q)} a $. It follows that $ (\m{D},Q) \in \Ob(\FA_{(\cdot)}^{Sg}) \circ \vDash_{(\cdot)}^{Sg}(\m{X}) $ and so $ \Ob(\FA^{Sg}_{(\cdot)}) \circ \vDash^{Sg}_{(\cdot)}(\m{X}) $ is stable under $ \mathbb{U}_{Sg} $. Therefore, $ \Ob(\FA^{Sg}_{(\cdot)}) \circ \vDash^{Sg}_{(\cdot)}(\m{X}) \supseteq \mathbb{U}_{Sg} \mathbb{S}\mathbb{P}(\m{X}) $.
		
		Now let $ (\m{D},Q) \in \Ob(\FA^{Sg}_{(\cdot)}) \circ \vDash^{Sg}_{(\cdot)}(\m{X}) $ and consider a sub-prop-morphism $ \iota \colon (\m{C}_{T},P_{T}) \to (\m{D},Q) $ such that $ \Sg(T) = Sg $. If $ a \in A_{Sg} \setminus \A(T) $, then the generic $ T $-model $ G $ in $ (\m{C}_{T},P_{T}) $ does not satisfy $ a $ and so $ T \nvDash_{(\m{C}_{T},P_{T})} a$. Since $ \vDash_{(\m{C}_{T},P_{T})} \; \supseteq \; \vDash_{(\m{D},Q)} \; \supseteq \; \vDash_{\m{X}} $ for each $ a \notin T $, there exists a $ T $-model $ F_{a} \colon (\m{C}_{Sg},P_{Sg}) \to (\m{C}_{a},P_{a}) $ such that $ (\m{C}_{a},P_{a}) \in \m{X} $ and $ G\db{a} \notin \ker F_{a} $. Then $ \ker(\epsilon_{T}) \leq \ker(\langle F_{a} \rangle_{a \notin T}) $ and so there exists a morphism $ \lambda \colon (\m{C}_{T},P_{T}) \to \prod_{a \notin T}(\m{C}_{a},P_{a}) $ such that
		\begin{equation*}
		\begin{tikzcd}
		(\m{C}_{Sg},P_{Sg}) \ar[swap]{dr}{\textstyle \epsilon_{T}} \ar{rr}{\textstyle \langle F_{a} \rangle_{a \notin T}} & & \prod_{a \notin T}(\m{C}_{a},P_{a}) \\
		& (\m{C}_{T},P_{T}) \ar[swap,dashed]{ru}{\textstyle \lambda } &
		\end{tikzcd}
		\end{equation*}
		commutes. Moreover, $ \Th(\epsilon_{T}) = \Th(\langle F_{a}\rangle_{a \notin T}) $ and so $ \lambda $ is a sub-prop-morphism. Thus $ (\m{C}_{T},P_{T}) \in \mathbb{S}\mathbb{P}(\m{X}) $ and so $ (\m{D},Q) \in \mathbb{U}_{Sg}\mathbb{S}\mathbb{P}(\m{X}) $. 
	\end{proof}
	Adjunction \ref{eqn: restricted adjunction.} is relevant when one wants to consider a logic over a fixed signature. However, in categorical logic, where internal logic is used, it makes sense to have a logic be independent of a particular signature. Thus we also provide a characterization of the closure $ \Ob(\FA_{\vDash_{(\cdot)}}) $ in Adjunction \ref{eqn: signature independent adjunction}. For $ \m{X} \subseteq \Ob(\FA) $, we define $ \mathbb{U}(\m{X}) $ by $ (\m{D},Q) \in \mathbb{U}(\m{X}) $, if each classifying sub-prop-category $ (\m{C}_{T},P_{T}) $ of $ (\m{D},Q) $ is in $ \m{X} $. That is, $ (\m{D},Q) \in \mathbb{U}(\m{X}) $ if and only if for each signature $ Sg $, $ (\m{D},Q) \in \mathbb{U}_{Sg}(\m{X}) $. For the following result, we must assume all prop-categories and signatures are small.
	\begin{thrm}\label{thm: Characterization of Logic Closure.}
		$ \Ob(\FA_{\vDash_{(\cdot)}}) = \mathbb{U} \mathbb{S} \mathbb{P}(\cdot) $.
	\end{thrm}
	\begin{proof}
		Clearly $ \Ob(\FA_{\vDash_{\m{X}}}) $ is stable under $ \mathbb{S} $ and $ \mathbb{P} $. Suppose $ (\m{D},Q) \in \mathbb{U}(\Ob(\FA_{\vDash_{\m{X}}})) $, and let $ T $ be the theory of $ (\m{D},Q) $. Then $ (\m{C}_{T},P_{T}) \in \FA_{\vDash_{\m{X}}} $ and since $ \vDash_{(\m{D},Q)} \, = \ \vDash_{(\m{C}_{T},P_{T})} $, $ (\m{D},Q) \in \Ob(\FA_{\vDash_{\m{X}}}) $.

		Now suppose $ (\m{D},Q) \in \Ob(\FA_{\vDash_{\m{X}}}) $. Then  $ \vDash_{(\m{D},Q)} \; \supseteq \; \vDash_{\m{X}} $. Following the proof of Theorem \ref{thm: Fixed signature adjunction}, one shows for each signature $ Sg $, $ (\m{D},Q) \in \mathbb{U}_{Sg}\mathbb{S}\mathbb{P}(\m{X}) $ and so $ (\m{D},Q) \in \mathbb{U}\mathbb{S}\mathbb{P}(\m{X}) $.
	\end{proof}

	\section{Conclusion}
	We have seen how the proof systems of a broad class of first-order logics naturally give rise to hyperdoctrine semantics, which can thus be seen as ``syntax in disguise". Historically, this was a derogatory expression for the extended algebraic semantics of propositional logics generated by the Lindenbaum-Tarski models of their theories. However, many properties of logics have been established by proving global/closure properties of their associated extended algebraic semantics and it is in the expanded space of hyperdoctrine semantics that we were able to establish the homomorphism theorems and provide an algebraic characterization of two logical closure operators. Moreover, it is in this extended hyperdoctrine semantics that Shirasu \cite{Shirasu1995} was able to lift an algebraic proof of the disjunction and existence properties of substructural propositional logics to their corresponding predicate logics. In our notation, the main argument in \cite{Shirasu1995} shows that $ \FA_{\m{L}} $, for the logics $ \m{L} $ in question, is stable under a ``gluing'' construction.
	
	Thus, in addition to the program of generalizing universal algebraic results to fibered algebras, we see lifting these theorems connecting properties of propositional logics to properties of their extended algebraic semantics as an important test for the hyperdoctrine approach. In particular, an interesting class of results to lift are those connecting interpolation results for logics with amalgamation results for their corresponding algebraic semantics, since, these results have been extensively studied \cite{Czelakowski1999} \cite{Metcalfe2014}, and include classical first-order logic through its algebraic semantics via cylindrical algebras.

	We have also seen that the adjoint conditions on quantifiers and equality, are not necessary to develop the $ 2 $-categorical hyperdoctrine semantics for first-order logics. This is significant because many natural notions of quantifiers fail to be adjoints. As the propositional theory of abstract algebraic logic does not fix a propositional signature nor does it demand certain connectives satisfy particular properties, we believe a first-order extension of abstract algebraic logic should do the same with respect to the quantifiers. It would be interesting to further explore modeling the vast landscape of nonclassical quantifiers. We hope the lifting of properties of propositional logics to first-order extensions via their hyperdoctrine semantics and subsequent analysis of what properties on the quantifiers are necessary will lead to a deeper theory and classification of nonclassical quantifiers.

	\section*{Acknowledgments}
	
	We would like to thank Adam P\v{r}enosil for his many insights and helpful suggestions on this paper. We also thank the anonymous referee for their feedback, which has greatly improved the exposition of the paper.

	\bibliographystyle{plain}
	\bibliography{Fibered_Universal_Algebra}
\end{document}